\documentclass[12pt]{amsart}
\usepackage[T1]{fontenc}
\usepackage[cp852]{inputenc}
\usepackage{amsmath}
\usepackage{amssymb}
\usepackage{moreverb}
\newtheorem{thm}{Theorem}[section]

\newtheorem{prop}[thm]{Proposition}
\newtheorem{rem}[thm]{Remark}
\newtheorem{example}[thm]{Example}
\usepackage{graphicx}

\usepackage{amsfonts}
\usepackage{amssymb}
\usepackage{eucal}

\theoremstyle{definition}
\numberwithin{equation}{section}
\begin{document}

\begin{center}
{\bf{A note on some generalized curvature tensor}}
\end{center}

\vspace{1mm}

\begin{center}
Ryszard Deszcz, Ma\l gorzata G\l ogowska, Marian Hotlo\'{s},
Miroslava Petrovi\'{c}-Torga\v{s}ev, and Georges Zafindratafa
\end{center}

\vspace{1mm}

\begin{center}
{\sl{Dedicated to the memory of Professor Krishan Lal Duggal (1929--2022)}}
\end{center}









\vspace{5mm}

\noindent
{\bf{Abstract.}} 
For any semi-Riemannian manifold $(M,g)$ we define 
some generalized curvature tensor $E$ as a linear combination 
of Kulkarni-Nomizu products formed by the metric tensor, 
the Ricci tensor and its square of given manifold. 
That tensor is closely related to quasi-Einstein spaces, 
Roter spaces and some Roter type spaces.

\section{Introduction}

Let $(M,g)$
be a semi-Riemannian manifold.
We denote by
$g$, $R$, $S$, $\kappa$ and $C$, 
the metric tensor, the Riemann-Christoffel curvature tensor, 
the Ricci tensor,
the scalar curvature and the Weyl conformal curvature tensor 
of $(M,g)$, respectively. Further, let $A \wedge B$ 
be the Kulkarni-Nomizu product of symmetric 
$(0,2)$-tensors $A$ and $B$.
Now we can define the $(0,2)$-tensors 
$S^{2}$ and $S^{3}$,
the $(0,4)$-tensors
$R \cdot S$, $C \cdot S$ and
$Q(A,B)$, 
and 
the $(0,6)$-tensors
$R \cdot R$, 
$R \cdot C$, 
$C \cdot R$, 
$C \cdot C$ 
and 
$Q(A,T)$, 
where $T$ is a generalized curvature tensor.
For precise definitions of the symbols used, 
we refer to Section 2
of this paper, as well as to
{\cite[Section 1] {2016_DGJZ}}, 
{\cite[Section 1] {2020_DGZ}},
{\cite[Chapter 6] {DHV2008}} and
{\cite[Sections 1 and 2] {DP-TVZ}}.

A semi-Riemannian manifold $(M,g)$, $\dim M = n \geq 2$, is said to be 
an {\sl Einstein manifold} \cite{Besse-1987},
or an {\sl Einstein space}, if at every point of $M$ 
its Ricci tensor $S$ is proportional to $g$, 
i.e., 
\begin{eqnarray}
S &=& \frac{\kappa}{n}\, g 
\label{2020.10.3.c}
\end{eqnarray}
on $M$,
assuming that the scalar curvature $\kappa$ is constant when $n = 2$. 
According to {\cite[p. 432] {Besse-1987}}
this condition is called the {\sl Einstein metric condition}.

Let $(M,g)$ be a semi-Riemannian manifold of dimension $\dim M = n \geq 3$.
We set
\begin{eqnarray}
E &=& g \wedge S^{2} + \frac{n-2}{2} \, S \wedge S - \kappa \, g \wedge S
+ \frac{\kappa ^{2} - \mathrm{tr}_{g} (S^{2})}{2(n-1)}
\, g \wedge g .
\label{2022.11.10.aaa}
\end{eqnarray}
It is easy to check
that the tensor $E$ is a generalized curvature tensor.
Further,
we define the subsets ${\mathcal{U}}_{R}$ and ${\mathcal U}_{S}$ of $M$ by 
${\mathcal{U}}_{R}  = 
\{x \in M\, | \, R - \frac{\kappa }{(n-1) n}\, G \neq 0\ \mbox {at}\ x \}$
and 
${\mathcal U}_{S} =  \{x \in M\, | \, 
S - \frac{\kappa }{n}\, g \neq 0\ \mbox {at}\ x \}$, respectively, 
where $G = \frac{1}{2}\,  g \wedge g$. 
If $n \geq 4$ then
we define the set ${\mathcal U}_{C} \subset M$ as the set of all points 
of $(M,g)$ at which which $C \neq 0$.
We note that if $n \geq 4$ then 
${\mathcal{U}}_{S} \cup {\mathcal{U}}_{C} = {\mathcal{U}}_{R}$ 
(see, e.g., \cite{DGHHY}).

An extension of the class of Einstein manifolds form
quasi-Einstein, $2$-quasi-Einstein and partially Einstein manifolds.
A semi-Riemannian manifold $(M,g)$, 
$\dim M = n \geq 3$, 
is said to be a {\sl quasi-Einstein manifold}, or a {\sl quasi-Einstein space},
if 
\begin{eqnarray}
\mathrm{rank}\, (S - \alpha\, g) &=& 1
\label{quasi02}
\end{eqnarray}
on ${\mathcal U}_{S} \subset M$, where $\alpha $ is some function 
on ${\mathcal U}_{S}$.
It is known that every non-Einstein warped product manifold 
$\overline{M} \times _{F} \widetilde{N}$
with a $1$-dimensional $(\overline{M}, \overline{g})$ base manifold and
a $2$-dimensional manifold $(\widetilde{N}, \widetilde{g})$
or an $(n-1)$-dimensional Einstein manifold
$(\widetilde{N}, \widetilde{g})$, $\dim M = n \geq 4$, 
and a warping function $F$,
is a  quasi-Einstein manifold (see, e.g., \cite{{Ch-DDGP}, {2016_DGJZ}}). 
A Riemannian manifold 
$(M,g)$, $\dim M = n \geq 3$, whose Ricci tensor has an eigenvalue
of multiplicity $n-1$ is a non-Einstein quasi-Einstein manifold 
(cf. {\cite[Introduction] {P47}}). 
We mention that quasi-Einstein manifolds arose during the study 
of exact solutions
of the Einstein field equations and the investigation 
on quasi-umbilical hypersurfaces 
of conformally flat spaces 
(see, e.g., \cite{{DGHSaw}, {2016_DGJZ}} and references therein). 
Quasi-Einstein hypersurfaces 
in semi-Riemannian spaces of constant curvature
were studied among others in
\cite{{DGHS}, {R102}, {P104}, {G6}} (see also \cite{DGHSaw} 
and references therein).
Quasi-Einstein manifolds satisfying some pseudosymmetry 
type curvature conditions  
were investigated recently in
\cite{{P119}, {Ch-DDGP}, {DGHHY}, {DGHZ01}, {DeHoJJKunSh}}. 
Quasi-Einstein hypersurfaces in conformally flat semi-Riemannian manifolds 
were studied in \cite{{Duggal-Sharma}, {Sharma}}.
In those papers quasi-Einstein hypersurfaces were called 
{\sl{pseudo-Einstein hypersurfaces}} (see also \cite{{KON}, {MAEDA}}).
Similarly, in \cite{{DVY-1994}, {GKZ-1991}} 
quasi-Einstein semi-Riemannian manifolds (hypersurfaces) 
were called {\sl{pseudo-Einstein manifolds (hypersurfaces)}}.

Let $(M,g)$, $\dim M = n \geq 3$, be a semi-Rieman\-nian manifold.
We note that (\ref{quasi02}) 
holds at a point $x \in {\mathcal U}_{S} \subset M$ 
if and only if 
$(S - \alpha\, g) \wedge (S - \alpha\, g) = 0$
at $x$, i.e.,  
\begin{eqnarray}
\frac{1}{2}\, S \wedge S   - \alpha\, g \wedge S 
+ \frac{\alpha ^{2}}{2}\, g \wedge g &=& 0    
\label{quasi03}
\end{eqnarray}
at $x$  (cf. {\cite[Proposition 2.1] {G6}}).
From (\ref{quasi03}), by a suitable contraction,
we get immediately 
\begin{eqnarray}
S^{2} 
&=& (\kappa - (n-2) \alpha )\, S + \alpha ( (n-1) \alpha - \kappa )\, g.
\label{quasi03quasi03} 
\end{eqnarray}

Using (\ref{2020.10.3.c}) we can easily check that the following equation
is satisfied on any Einstein manifold $(M,g)$
\begin{eqnarray}
g \wedge S^{2} + \frac{n-2}{2} \, S \wedge S - \kappa \, g \wedge S
+ \frac{\kappa ^{2} - \mathrm{tr}_{g} (S^{2})}{2(n-1)} \, g \wedge g 
&=& 0 ,
\label{2022.11.20.aa}
\end{eqnarray}
i.e., $E = 0$ on $M$, where the tensor $E$ 
is defined by  (\ref{2022.11.10.aaa}).
Moreover, 
as it was stated in {\cite[Lemma 2.1] {DGHSaw-2022}},
(\ref{2022.11.20.aa})  is satisfied on every quasi-Einstein manifold 
$(M,g)$, $n \geq 3$. The converse statement also is true. Precisely,
from Proposition 2.1 it follows that 
if  $(M,g)$, $n \geq 4$, is a semi-Riemannian manifold satisfying  
(\ref{2022.11.20.aa}) on ${\mathcal U}_{S} \subset M$ then 
a condition of the form (\ref{quasi02}) holds on ${\mathcal U}_{S}$,
where $\alpha $ is some function on this set.
In Section 2 we also present another result 
related to the tensor $E$ (see Proposition 2.2). 
Namely, we prove that if a generalized curvature tensor $T$
is a linear combination of the tensors 
$R$, $S\wedge S$, $g \wedge S$, $g \wedge S^{2}$, 
and $g \wedge g$ then the Weyl tensor of $T$ is a linear combination
of the tensors $C$ and $E$.  
The tensor $E$ is determined by some Kulkarni-Nomizu products formed 
by $g$, $S$ and $S^{2}$,
i.e.,  $E$ is defined by (\ref{2022.11.10.aaa}). 
In the same way, we can define the $(0,4)$-tensor $E(A)$ corresponding 
to a symmetric  $(0,2)$-tensor $A$ 
\begin{eqnarray}
E(A) &=& g \wedge A^{2} + \frac{n-2}{2} \, A \wedge A 
- \mathrm{tr}_{g} ( A)  \, g \wedge A
+ \frac{(\mathrm{tr}_{g} ( A) )^{2} - \mathrm{tr}_{g} (A^{2})}{2(n-1)}
\, g \wedge g .
\label{2022.12.33.aaa}
\end{eqnarray}

The semi-Riemannian manifold $(M,g)$, $\dim M = n \geq 3$, will be called
a {\sl{partially Einstein manifold}}, or
a {\sl{partially Einstein space}}
(cf. {\cite[Foreword] {CHEN-2017}}, {\cite[p. 20] {V2}}),  
if at every point $x \in {\mathcal{U}}_{S} \subset M$ 
its Ricci operator ${\mathcal{S}}$ satisfies
${\mathcal{S}}^{2} = \lambda {\mathcal{S}} + \mu  I\!d_{x}$,
or equivalently, 
\begin{eqnarray}
S^{2} &=&  \lambda \, S + \mu \, g ,
\label{partiallyEinstein}
\end{eqnarray}
where $\lambda , \mu \in {\mathbb{R}}$  
and $I\!d_{x}$ is the identity transformation of $T_{x} M$.
Evidently, 
(\ref{quasi03quasi03})
is a special case of (\ref{partiallyEinstein}). 
Thus every quasi-Einstein manifold is a partially Einstein manifold.
The converse statement is not true. 
Contracting (\ref{partiallyEinstein}) we get
$ \mathrm{tr}_{g} (S^{2}) = \lambda \, \kappa + n\, \mu$. 
This together with 
(\ref{partiallyEinstein}) yields (cf. {\cite[Section 5] {2021-DGH}})
\begin{eqnarray*} 
S^{2} - \frac{ \mathrm{tr}_{g} (S^{2})}{n} \, g 
&=&  \lambda \left( S - \frac{\kappa }{n} \, g \right) .
\end{eqnarray*}
In particular, a Riemannian manifold $(M,g)$, $\dim M = n \geq 3$, 
is a partially Einstein space if at every point 
$x \in {\mathcal{U}}_{S} \subset M$ 
its Ricci operator ${\mathcal{S}}$
has exactly two distinct eigenvalues
$\kappa _{1}$ and $\kappa _{2}$ 
with multiplicities $p$ and $n-p$, respectively, where
$1 \leq p \leq n-1$.
Evidently, if $p = 1$, or $p = n-1$, 
then $(M,g)$ is a quasi-Einstein manifold. 

In Section 3 we present definitions of some classes of semi-Riemannian manifolds  
determined by curvature conditions of pseudosymmetry type. 
Investigations of semi-Riemannian manifolds satisfying some particular
curvature conditions of pseudosymmetry type lead to 
Roter spaces (see Propositon 4.1).
Roter spaces form an important subclass 
of the class of non-conformally flat and non-quasi-Einstein
partially Einstein manifolds of dimension $\geq 4$.
Namely, a non-quasi-Einstein and non-conformally flat 
semi-Riemannian manifold $(M,g)$, $\dim M = n \geq 4$, 
satisfying on 
${\mathcal U}_{S} \cap {\mathcal U}_{C} \subset M$ 
the following equation
\begin{eqnarray}
R &=& \frac{\phi}{2}\, S\wedge S + \mu\, g\wedge S 
+ \frac{\eta}{2}\, g \wedge g ,
\label{eq:h7a}
\end{eqnarray}
where 
$\phi$, $\mu $ and $\eta $ are some functions on this set,
is called a {\sl Roter type manifold}, or a {\sl Roter manifold}, 
or a {\sl Roter space} 
(see, e.g., {\cite[Section 15.5] {Chen-2021}},
\cite{{P106}, {2016_DGJZ}, {DGP-TV02}, {DHV2008}}).
Equation (\ref{eq:h7a}) is called a {\sl{Roter equation}}
(see, e.g., {\cite[Section 1] {DGHSaw-2022}}). 
In Section 4 we present results on such manifolds. 
For instance,  every Roter space $(M,g)$, $\dim M = n \geq 4$, 
satisfies among others the following pseudosymmetry type curvature condition 
on ${\mathcal{U}}_{S} \cap {\mathcal{U}}_{C} \subset M$
(see Theorem 4.2) 
\begin{eqnarray}
C \cdot R - R \cdot C = Q(S,C) - \frac{\kappa}{ n-1 }\, Q(g,C) .
\label{Roterformula}
\end{eqnarray}

Let $(M,g)$, $\dim M = n \geq 4$, be a non-partially-Einstein 
and non-conformally flat semi-Rieman\-nian manifold.
If its Riemann-Christoffel curvature $R$
is at every point of  
${\mathcal U}_{S} \cap {\mathcal U}_{C} \subset M$ 
a linear combination of the Kulkarni-Nomizu products 
formed by the tensors
$S^{0} = g$ and $S^{1} = S, \ldots , S^{p-1}, S^{p}$, 
where $p$ is some natural 
number $\geq 2$, then $(M,g)$
is called 
a {\sl generalized Roter type manifold},
or a {\sl generalized Roter manifold}, 
or a {\sl generalized Roter type space}, or 
a {\sl generalized Roter space}. For instance, when $p = 2$,
we have
\begin{eqnarray} 
R &=&  \frac{\phi _{2} }{2}\, S^{2} \wedge S^{2} 
+ \phi_{1}\, S \wedge S^{2}  
+ \frac{\phi}{2}\, S \wedge S
+ \mu _{1}\, g \wedge S^{2}
+ \mu \, g \wedge S
+ \frac{ \eta }{2} \, g \wedge g ,
\label{B001simply}
\end{eqnarray}
where 
$\phi $, 
$\phi _{1}$,
$\phi _{2}$,
$\mu_{1}$, $\mu$ and 
$\eta$ 
are functions on
${\mathcal U}_{S} \cap {\mathcal U}_{C}$.
Because $(M,g)$ is a non-partially Einstein manifold,
at least one of the functions 
$\mu _{1}$, $\phi _{1}$ and $\phi _{2}$ is a non-zero function.
Equation (\ref{B001simply}) is called a 
{\sl{Roter type equation}}  
(see, e.g., {\cite[Section 1] {DGHSaw-2022}}).
We refer to 
\cite{{DGHSaw-2022}, {2013_DGJP-TZ}, {2016_DGJZ}, 
{DGP-TV02}, {DeHoJJKunSh}, 
{S3}, {Saw114}, {SDHJK}, {2016_SK}, {2019_SK}}
for results on manifolds (hypersurfaces) satisfying (\ref{B001simply}).

As it was stated in {\cite[Lemma 2.2] {DGHSaw-2022}} (see Proposition 4.4),
if $(M,g)$, $\dim M = n \geq 4$, is a Roter space satisfying
(\ref{eq:h7a}) 
on ${\mathcal U}_{S} \cap  {\mathcal U}_{C} \subset M$ 
then on this set
\begin{eqnarray}
C &=& \frac{\phi}{n-2} 
\left(g \wedge S^{2} + \frac{n-2}{2} \, S \wedge S - \kappa \, g \wedge S
+ \frac{\kappa ^{2} - \mathrm{tr}_{g} (S^{2})}{2(n-1)}
\, g \wedge g \right) .
\label{2022.11.22.aa}
\end{eqnarray}

In Section 5 we recall results on some warped product manifolds 
with $2$-dimensional base manifold obtained in \cite{DGHSaw-2022}.

In Section 6 we state that on every essentially conformally symmetric manifold
the following equation is satisfied
\begin{eqnarray}
F\, C 
&=& \frac{n-2}{2 (n -2)} \, S \wedge S\nonumber\\
&=& \frac{1}{n-2} 
\left(  
g \wedge S^{2} + \frac{n-2}{2} \, S \wedge S - \kappa \, g \wedge S
+ \frac{\kappa ^{2} - \mathrm{tr}_{g} (S^{2})}{2(n-1)} \, g \wedge g 
\right) .
\label{2023.03.03.a}
\end{eqnarray}

In Section 7 we recall some known results on hypersurfaces $M$,
$\dim M \geq 4$, isometrically immersed in a conformally flat spaces. 
In particular, 
we mention that 
at every point of $M$ its Weyl conformal curvature tensor $C$  
and the $(0,4)$-tensor $E(H)$, formed for the second fundamental tensor $H$
of $M$, are linearly dependent (see (\ref{2022.08.03.dd})).

In the last section we investigate non-Einstein and non-conformally flat 
hypersurfaces $M$, $\dim M \geq 4$, isometrically immersed 
in semi-Riemannian spaces of constant curvature
satisfying some curvature conditions of pseudosymmetry type. 
Under some additional assumptions imposed on the second fundamental tensor 
$H$ of $M$ we obtain equations involved with the tensor $E$.

\section{Preliminaries.}

Throughout this paper, all manifolds are assumed 
to be connected paracompact
mani\-folds of class $C^{\infty }$. Let $(M,g)$, $\dim M = n \geq 3$,
be a semi-Riemannian manifold, and let $\nabla$ 
be its Levi-Civita connection and $\Xi (M)$ the Lie
algebra of vector fields on $M$. We define on $M$ the endomorphisms 
$X \wedge _{A} Y$ and ${\mathcal{R}}(X,Y)$ of $\Xi (M)$ by
$(X \wedge _{A} Y)Z = A(Y,Z)X - A(X,Z)Y$ and
\begin{eqnarray*}
{\mathcal R}(X,Y)Z 
&=& 
\nabla _X \nabla _Y Z - \nabla _Y \nabla _X Z - \nabla _{[X,Y]}Z ,
\end{eqnarray*}
respectively, where 
$X, Y, Z \in \Xi (M)$
and 
$A$ is a symmetric $(0,2)$-tensor on $M$. 
The Ricci tensor $S$, the Ricci operator ${\mathcal{S}}$ 
and the scalar curvature
$\kappa $ of $(M,g)$ are defined by 
\begin{eqnarray*}
S(X,Y)\ =\ \mathrm{tr} \{ Z \rightarrow {\mathcal{R}}(Z,X)Y \} ,\ \  
g({\mathcal S}X,Y)\ =\ S(X,Y) ,\ \  
\kappa \ =\ \mathrm{tr}\, {\mathcal{S}},
\end{eqnarray*}
respectively. 
The endomorphism ${\mathcal{C}}(X,Y)$ is defined by
\begin{eqnarray*}
{\mathcal C}(X,Y)Z  &=& {\mathcal R}(X,Y)Z 
- \frac{1}{n-2}(X \wedge _{g} {\mathcal S}Y + {\mathcal S}X \wedge _{g} Y
- \frac{\kappa}{n-1}X \wedge _{g} Y)Z .
\end{eqnarray*}
Now the $(0,4)$-tensor $G$, 
the Riemann-Christoffel curvature tensor $R$ and
the Weyl conformal curvature tensor $C$ of $(M,g)$ are defined by
$G(X_1,X_2,X_3,X_4) = g((X_1 \wedge _{g} X_2)X_3,X_4)$ and
\begin{eqnarray*}
R(X_1,X_2,X_3,X_4) &=& g({\mathcal R}(X_1,X_2)X_3,X_4) ,\\
C(X_1,X_2,X_3,X_4) &=& g({\mathcal C}(X_1,X_2)X_3,X_4) ,
\end{eqnarray*}
respectively, where $X_1,X_2,\ldots \in \Xi (M)$.  
For a symmetric $(0,2)$-tensor $A$ we denote by 
${\mathcal{A}}$ the endomorphism related to $A$ by 
$g({\mathcal{A}}X,Y) = A(X,Y)$.
The $(0,2)$-tensors
$A^{p}$, $p = 2, 3, \ldots $, are defined by 
$A^{p}(X,Y) = A^{p-1} ({\mathcal{A}}X, Y)$,
assuming that $A^{1} = A$. In this way, for $A = S$ and 
${\mathcal{A}} = {\mathcal S}$
we get the tensors $S^{p}$, 
$p = 2, 3, \ldots $, assuming that $S^{1} = S$.

Let ${\mathcal B}$ be a tensor field sending any $X, Y \in \Xi (M)$ 
to a skew-symmetric endomorphism 
${\mathcal B}(X,Y)$, 
and let $B$ be
the $(0,4)$-tensor associated with ${\mathcal B}$ by
\begin{eqnarray}
B(X_1,X_2,X_3,X_4) &=& 
g({\mathcal B}(X_1,X_2)X_3,X_4) .
\label{DS5}
\end{eqnarray}
The tensor $B$ is said to be a {\sl{generalized curvature tensor}}  if the
following two conditions are fulfilled:
$B(X_1,X_2,X_3,X_4) = B(X_3,X_4,X_1,X_2)$ and
\begin{eqnarray*}
B(X_1,X_2,X_3,X_4) + B(X_2,X_3,X_1,X_4) + B(X_3,X_1,X_2,X_4) &=& 0 . 
\end{eqnarray*}
For ${\mathcal B}$ as above, let $B$ be again defined by (\ref{DS5}). 
We extend the endomorphism ${\mathcal B}(X,Y)$
to a derivation ${\mathcal B}(X,Y) \cdot \, $ of the algebra 
of tensor fields on $M$,
assuming that it commutes with contractions and 
${\mathcal B}(X,Y) \cdot \, f  = 0$ for any smooth function $f$ on $M$. 
Now for a $(0,k)$-tensor field $T$,
$k \geq 1$, we can define the $(0,k+2)$-tensor $B \cdot T$ by
\begin{eqnarray*}
& & (B \cdot T)(X_1,\ldots ,X_k,X,Y) \ =\ 
({\mathcal B}(X,Y) \cdot T)(X_1,\ldots ,X_k)\\  
&=& - T({\mathcal{B}}(X,Y)X_1,X_2,\ldots ,X_k)
- \cdots - T(X_1,\ldots ,X_{k-1},{\mathcal{B}}(X,Y)X_k) .
\end{eqnarray*}
If $A$ is a symmetric $(0,2)$-tensor then we define the
$(0,k+2)$-tensor $Q(A,T)$ by
\begin{eqnarray*}
& & Q(A,T)(X_1, \ldots , X_k, X,Y) \ =\
(X \wedge _{A} Y \cdot T)(X_1,\ldots ,X_k)\\  
&=&- T((X \wedge _A Y)X_1,X_2,\ldots ,X_k) 
- \cdots - T(X_1,\ldots ,X_{k-1},(X \wedge _A Y)X_k) .
\end{eqnarray*}
In this manner we obtain the $(0,6)$-tensors 
$B \cdot B$ and $Q(A,B)$.

Substituting in the above formulas 
${\mathcal{B}} = {\mathcal{R}}$ or ${\mathcal{B}} = {\mathcal{C}}$, 
$T=R$ or $T=C$ or $T=S$, $A=g$ or $A=S$ 
we get the tensors $R\cdot R$, $R\cdot C$, $C\cdot R$, 
$C\cdot C$, $R\cdot S$, 
$Q(g,R)$, $Q(S,R)$, $Q(g,C)$,  $Q(S,C)$, 
and $Q(g,S)$, $Q(g,S^{2})$,  $Q(S,S^{2})$.

For a symmetric $(0,2)$-tensor $A$ and a $(0,k)$-tensor $T$, $k \geq 2$, we
define their {\sl{Kulkarni-Nomizu tensor}} $A \wedge T$ by 
(see, e.g., {\cite[Section 2] {DGHHY}}) 
\begin{eqnarray*}
& &(A \wedge T )(X_{1}, X_{2}, X_{3}, X_{4}; Y_{3}, \ldots , Y_{k})\\
&=&
A(X_{1},X_{4}) T(X_{2},X_{3}, Y_{3}, \ldots , Y_{k})
+ A(X_{2},X_{3}) T(X_{1},X_{4}, Y_{3}, \ldots , Y_{k} )\\
& &
- A(X_{1},X_{3}) T(X_{2},X_{4}, Y_{3}, \ldots , Y_{k})
- A(X_{2},X_{4}) T(X_{1},X_{3}, Y_{3}, \ldots , Y_{k}) .
\end{eqnarray*}
It is obvious that the following tensors 
are generalized curvature tensors: $R$, $C$ and 
$A \wedge B$, where $A$ and $B = T$ are symmetric $(0,2)$-tensors. 
We have 
\begin{eqnarray}
C &=& R - \frac{1}{n-2}\, g \wedge S + \frac{\kappa }{(n-2) (n-1)} \, G , 
\label{Weyl}\\
G &=& \frac{1}{2}\, g \wedge g ,
\label{2020.12.08.a}
\end{eqnarray}
and (see, e.g., {\cite[Lemma 2.2(i)] {DGHHY}})
\begin{eqnarray}
&(a)&  
\ \
Q(A, A \wedge B) \ =\ - \frac{1}{2}\, Q(B, A \wedge A ) ,\nonumber\\
&(b)&
\ \
A \wedge Q(A,B) \ =\ - \frac{1}{2}\, Q(B, A \wedge A ).
\label{DS7}
\end{eqnarray}
By an application of (\ref{DS7})(a) we obtain on $M$ the identities
\begin{eqnarray}
Q(g, g \wedge S) &=& - Q(S,G) ,
\ \ Q(S, g \wedge S)\ = \ - \frac{1}{2}\, Q(g, S \wedge S) .  
\label{dghz01}
\end{eqnarray}
Further, by making use of (\ref{Weyl}), (\ref{2020.12.08.a})
and (\ref{dghz01}), we get immediately
\begin{eqnarray}
Q(g,C) &=& Q(g,R) - \frac{1}{n-2}\, Q(g,g \wedge S) 
+ 
 \frac{\kappa }{(n-2) (n-1)} \, Q(g,G)  \nonumber\\
&=& Q(g,R) - \frac{1}{n-2}\, Q(g,g \wedge S) 
\ = \ Q(g,R) + \frac{1}{n-2} \, Q(S, \frac{1}{2}\, g \wedge g) ,
\label{abcd07zz},\\ 
Q(S,C) &=& Q(S,R) - \frac{1}{n-2}\, Q(S,g \wedge S) 
+ 
 \frac{\kappa }{(n-2) (n-1)} \, Q(S,G) \nonumber\\ 
&=&  
Q(S,R) + \frac{1}{2 (n-2)}\, Q(g,S \wedge S) 
- \frac{\kappa }{(n-2) (n-1)} \, Q(g,g \wedge S) .
\label{abcd07}
\end{eqnarray}

From (\ref{DS7}) (a) it follows immmediately that $Q(g, g \wedge g) = 0$.
Thus we have
\begin{eqnarray}
Q(g, E) &=& 
Q(g, g \wedge S^{2} + \frac{n-2}{2} \, S \wedge S - \kappa \, g \wedge S ) ,
\label{2022.11.10.bbb}
\end{eqnarray}
where the tensor $E$ is defined by  (\ref{2022.11.10.aaa}).

Let $A_{1}$, $A_{2}$ and $B$ be
symmetric $(0,2)$-tensors.
We have (see, e.g., 
{\cite[Lemma 2.1(i)] {Ch-DDGP}} and references therein) 
\begin{eqnarray}
A_{1} \wedge Q(A_{2},B) + A_{2} \wedge Q(A_{1},B) 
+ Q(B, A_{1} \wedge A_{2}) 
&=& 0 .
\label{2021.11.30.a4}
\end{eqnarray} 
From (\ref{2021.11.30.a4}) we get easily
(see also {\cite[Lemma 2.2(iii)] {DGHHY}}
and references therein)
\begin{eqnarray*}
 Q(B, A_{1} \wedge A_{2}) 
+ Q(A_{1}, A_{2} \wedge B)
+ Q(A_{2}, B \wedge A_{1}) &=& 0 .
\end{eqnarray*}

Let $A$ be a symmetric $(0,2)$-tensor and 
$T$ a $(0,k)$-tensor, $k = 2, 3, \ldots $.  
The tensor $Q(A,T)$ is called the {\sl{Tachibana tensor}} of $A$ and $T$, 
or the Tachibana tensor for short (see, e.g., \cite{DGPSS}). 
Using the tensors $g$, $R$ and $S$ we can define the following 
$(0,6)$-Tachibana tensors: 
$Q(S,R)$, $Q(g,R)$, $Q(g,g \wedge S)$ and $Q(S, g \wedge S)$. 
We can check, by making use of (\ref{DS7})(a) and (\ref{dghz01}), 
that other $(0,6)$-Tachibana tensors  
constructed from $g$, $R$ and $S$ may be expressed 
by the four Tachibana tensors mentioned above 
or vanish identically on $M$.

Let $T$ be a generalized curvature tensor on a semi-Riemannian manifold $(M,g)$, 
$\dim M = n \geq 4$. We denote by
$\mathrm{Ric} (T)$, $\kappa (T)$ and $\mathrm{Weyl} (T)$ the Ricci tensor, 
the scalar curvature and the Weyl tensor of the tensor $T$, respectively.
We refer to {\cite[Section 2] {DGHHY}},  {\cite[Section 3] {2021-DGH}} 
or {\cite[Section 3] {DGHZ01}} 
for definitions of the considered tensors.
In particular, we have 
\begin{eqnarray}
\mathrm{Weyl} (T) &=& T - \frac{1}{n-2} \, g \wedge \mathrm{Ric} (T) 
+  \frac{\kappa (T) }{2 (n-2)(n-1)}\, g \wedge g .
\label{2023.03.06a}
\end{eqnarray}

Let $A$ be a symmetric $(0,2)$-tensor on a semi-Riemannian manifold
$(M,g)$, $\dim M = n \geq 3$. Let $E(A)$ be the tensor 
defined by (\ref{2022.12.33.aaa}). It is easy to check that 
$\mathrm{Ric} (E(A))$ is a zero tensor. Therefore, we also have
$\kappa (E(A)) = 0$. Any generalized curvature tensor $T$
defined on a $3$-dimensional semi-Riemannian manifold $(M,g)$ is expressed by
$T = g \wedge \mathrm{Ric} (T) - (\kappa (T)/4) g \wedge g$
{\cite[p. 48] {Eisenhart}} (see also {\cite[Lemma 2 (ii)] {D-1991}}).
Thus we see that the tensor $T = E(A)$ on any $3$-dimensional 
semi-Riemannian manifold $(M,g)$ is a zero tensor. In particular, 
on any $3$-dimensional semi-Riemannian manifold $(M,g)$ we have
$E = 0$.

Let $A$ be a symmetric $(0,2)$-tensor on a semi-Riemannian manifold
$(M,g)$, $\dim M = n \geq 3$. We denote by ${\mathcal U}_{A}$ the set
of points of $M$ at which $A \neq \frac{\mathrm{tr}_{g} (A)}{n}\, g$.
\begin{prop}
Let $A$ be a symmetric $(0,2)$-tensor on a semi-Riemannian manifold
$(M,g)$, $\dim M = n \geq 4$.
(i) (cf. {\cite[Lemma 2.1] {DGHSaw-2022}}) 
If the following condition is satisfied on ${\mathcal U}_{A}  \subset M$
\begin{eqnarray}
\mathrm{rank} ( A - \alpha \, g ) &=& 1 
\label{2022.12.22.aa} 
\end{eqnarray}
then
\begin{eqnarray}
g \wedge A^{2} + \frac{n-2}{2}\, A \wedge A 
- \mathrm{tr}_{g} (A) \, g \wedge A 
+ \frac{(\mathrm{tr}_{g} (A))^{2} 
- \mathrm{tr}_{g} (A^{2})}{2 (n-1)} \, g \wedge g 
&=& 0 
\label{2022.12.22.bb} 
\end{eqnarray}
and
\begin{eqnarray}
A^{2} - \frac{ \mathrm{tr}_{g} (A^{2})}{n} &=&
( \mathrm{tr}_{g} (A) - (n - 2) \alpha ) 
\left( A - \frac{\mathrm{tr}_{g} (A) }{n} \, g \right)  
\label{2022.12.20.kk} 
\end{eqnarray}
on ${\mathcal U}_{A}$,
where $\alpha$ is some function on ${\mathcal U}_{A}$. 
\newline
(ii) If (\ref{2022.12.22.bb}) is satisfied 
on ${\mathcal U}_{A} \subset M$ then 
\begin{eqnarray}
A^{2} - \frac{\mathrm{tr}_{g} (A^{2})}{n}\, g 
&=& 
\rho \left( A - \frac{\mathrm{tr}_{g} (A)}{n}\, g \right) 
\label{2022.12.22.ii} 
\end{eqnarray}
and
\begin{eqnarray}
\left(
A - \frac{ \mathrm{tr}_{g} (A) - \rho }{n-2}\, g
\right)
\wedge
\left(
A - \frac{ \mathrm{tr}_{g} (A) - \rho }{n-2}\, g
\right) &=& 0
\label{2022.12.22.jj} 
\end{eqnarray}
on ${\mathcal U}_{A}$,
where $\rho$ is some function on ${\mathcal U}_{A}$.
\end{prop}
{\bf{Proof.}} 
(i) (cf. the proof of {\cite[Lemma 2.1] {DGHSaw-2022}})
Equation (\ref{2022.12.22.aa}) yields {\cite[Proposition 2.2] {G6}}
\begin{eqnarray}
\frac{1}{2}\, A \wedge A &=& \alpha \, g \wedge A 
- \frac{\alpha ^{2} }{2}\, g \wedge g .
\label{2022.12.20.ll} 
\end{eqnarray}
This, by suitable contractions yields
\begin{eqnarray}
A^{2} - \mathrm{tr}_{g} (A)\, A &=&
- (n-2) \alpha \, A - \alpha \mathrm{tr}_{g} (A)\, g 
+ (n-1) \alpha ^{2}\, g ,
\label{2022.12.20.yy}\\ 
\mathrm{tr}_{g} (A^{2}) - (\mathrm{tr}_{g} (A))^{2}
&=& - 2(n-1) \alpha \mathrm{tr}_{g} (A) + n (n-1) \alpha ^{2} .
\label{2022.12.20.zz} 
\end{eqnarray}
Now using (\ref{2022.12.20.ll}), (\ref{2022.12.20.yy})
and (\ref{2022.12.20.zz}) we can easily check that 
(\ref{2022.12.22.bb}) and (\ref{2022.12.20.kk}) 
hold on ${\mathcal U}_{A}$.
\newline
(ii) (cf. the proof of {\cite[Lemma 3.4] {D-1990}})
In the local coordinates (\ref{2022.12.22.bb}) reads 
\begin{eqnarray}
& &
g_{hk}A^{2}_{ij} + g_{ij}A^{2}_{hk} - g_{hj}A^{2}_{ik} - g_{ik}A^{2}_{hj}
+
(n-2)\, (A_{hk}A_{ij} - A_{hj}A_{ik}) 
\nonumber\\
& &
- \mathrm{tr}_{g} (A) \,
(g_{hk}A_{ij} + g_{ij}A_{hk} - g_{hj}A_{ik} - g_{ik}A_{hj})\nonumber\\
& &
+ \frac{(\mathrm{tr}_{g} (A))^{2} - \mathrm{tr}_{g} (A^{2})}{n-1}
\, (g_{hk}g_{ij} - g_{hj}g_{ik}) \ = \ 0 .
\label{2022.12.22.cc} 
\end{eqnarray}
Contracting (\ref{2022.12.22.cc})  with 
$A^{ij} = A_{rs}g^{ri}g^{sj}$ and $A^{k}_{\ l} = A_{rl}g^{rk}$ we find
\begin{eqnarray}
A^{3} &=& \frac{3 \mathrm{tr}_{g} (A)}{n}\, A^{2}
+ \left(
\frac{(n^{2} - 3n + 3) \mathrm{tr}_{g} (A^{2})}{ (n-1) n } 
- \frac{(\mathrm{tr}_{g} (A) )^{2}}{n-1} 
\right) A\nonumber\\
& & 
+ \left(
\frac{  (\mathrm{tr}_{g} (A))^{3} }{(n-1) n}
- \frac{ \mathrm{tr}_{g} (A) \, \mathrm{tr}_{g} (A^{2})}{n-1}
+ \frac{ \mathrm{tr}_{g} (A^{3})}{n}
\right) g ,
\label{2022.12.22.dd} 
\end{eqnarray}
\begin{eqnarray}
& &
A_{hl}A^{2}_{ij} - A_{il}A^{2}_{hj}
+ g_{ij}A^{3}_{hl} - g_{hj}A^{3}_{il} 
+ (n-2)\, (A_{ij}A^{2}_{hl} - A_{hj}A^{2}_{il}) 
\nonumber\\
& &
- \mathrm{tr}_{g} (A) \,
(A_{hl}A_{ij} - A_{il}A_{hj} 
+ g_{ij}A^{2}_{hl} - g_{hj}A^{2}_{il} )\nonumber\\
& &
+ \frac{(\mathrm{tr}_{g} (A))^{2} - \mathrm{tr}_{g} (A^{2})}{n-1}
\, (g_{ij}A_{hl} - g_{hj}A_{il}) \ = \ 0 ,
\label{2022.12.22.ee} 
\end{eqnarray}
respectively. 
From (\ref{2022.12.22.ee}), by symmetrization in $l,j$, we obtain
\begin{eqnarray}
& &
Q(g, A^{3}) + (n-3) Q(A,A^{2}) - \mathrm{tr}_{g} (A) \, Q(g, A^{2})\nonumber\\
& &
+
\frac{(\mathrm{tr}_{g} (A))^{2} 
- \mathrm{tr}_{g} (A^{2})}{n-1} \, Q(g,A) \ =\ 0 .
\label{2022.12.22.ff} 
\end{eqnarray} 
Applying (\ref{2022.12.22.dd}) into (\ref{2022.12.22.ff}) we get
\begin{eqnarray*}
(n-3) \left(
Q(A,A^{2}) - \frac{\mathrm{tr}_{g} (A)}{n}\, Q(g, A^{2})
+ \frac{\mathrm{tr}_{g} (A^{2})}{n}\, Q(g, A) \right)
&=& 0 ,
\end{eqnarray*} 
which yields
\begin{eqnarray*}
Q\left( A - \frac{\mathrm{tr}_{g} (A)}{n}\, g , 
A^{2} - \frac{\mathrm{tr}_{g} (A^{2})}{n}\, g \right) &=& 0 .
\end{eqnarray*} 
From this, in view of {\cite[Lemma 2.4 (ii)] {DV-1991}},
it follows that (\ref{2022.12.22.ii}) holds on ${\mathcal U}_{A}$.
Now (\ref{2022.12.22.bb}) and (\ref{2022.12.22.ii}),
by an application of {\cite[Lemma 3.1] {G2}}, 
lead to (\ref{2022.12.22.jj}), completing the proof of (ii).


\begin{prop}
Let $T$ be a generalized curvature tensor on a semi-Riemannian manifold $(M,g)$, 
$\dim M = n \geq 4$.
If the following condition is satisfied at a point $x \in M$ 
\begin{eqnarray}
T &=& \alpha _{1} \, R + \frac{\alpha _{2}}{2} \, S \wedge S 
+ \alpha _{3}\, g \wedge S  +  \alpha _{4}\, g \wedge S^{2} 
+  \frac{ \alpha _{5} }{2} \, g \wedge g
\label{2023.03.06b}
\end{eqnarray}
then 
\begin{eqnarray}
\mathrm{Weyl} (T) &=&  \alpha _{1} \, C +   \frac{\alpha _{2}}{n-2} \, E 
\label{2023.03.06c}
\end{eqnarray}
at this point, where the tensor $E$ is defined by (\ref{2022.11.10.aaa}) 
and $\alpha_{1}, \ldots , \alpha_{5} \in \mathbb{R}$.
\end{prop}
{\bf{Proof.}} From (\ref{2023.03.06b}), by a suitable contraction, 
we get immediately 
\begin{eqnarray}
\mathrm{Ric} (T) 
&=& ( \alpha _{1} + \alpha _{2} \kappa + (n-2)\, \alpha _{3} )\, S
 +  ( (n-2) \alpha _{4} - \alpha _{2} ) \, S^{2}
+ \alpha _{6} \, g  \, ,
\label{2023.03.06d}
\end{eqnarray}
where $\alpha _{6}$ is some real number. Now using
(\ref{2022.11.10.aaa}),  (\ref{Weyl}), (\ref{2020.12.08.a}),
(\ref{2023.03.06a}), (\ref{2023.03.06b}) and (\ref{2023.03.06d}) we get
\begin{eqnarray*}
& & 
\mathrm{Weyl} (T) \ =\ T - \frac{1}{n-2} \, g \wedge \mathrm{Ric} (T) 
+  \frac{\kappa (T) }{2 (n-2)(n-1)}\, g \wedge g \nonumber\\
&=&   \alpha _{1} \, R + \frac{\alpha _{2}}{2} \, S \wedge S 
+ \alpha _{3}\, g \wedge S  +  \alpha _{4}\, g \wedge S^{2}
+ \frac{\alpha _{7}}{2} \, g \wedge g  \nonumber\\
& &
- \frac{1}{n-2}  ( \alpha _{1} + \alpha _{2} \kappa 
+ (n-2)\, \alpha _{3} )\, g \wedge S
-  \frac{1}{n-2}  ( (n-2) \alpha _{4} - \alpha _{2} ) 
\, g \wedge S^{2} \nonumber\\
&=&
 \alpha _{1} \, R +  \frac{\alpha _{2}}{2} \, S \wedge S 
- \frac{ \alpha _{1} + \alpha _{2} \kappa }{n-2} \, g \wedge S
+  \frac{ \alpha _{2}}{n-2}  \, g \wedge S^{2}
+ \frac{\alpha _{7}}{2} \, g \wedge g \nonumber\\
&=&
 \alpha _{1} \, ( R - \frac{ 1 }{n-2} \, g \wedge S )
+ 
\frac{\alpha _{2}}{n-2}\, (  g \wedge S^{2} + \frac{n-2}{2}\, S \wedge S 
- \kappa \, g \wedge S)
+ \frac{\alpha _{7}}{2} \, g \wedge g \nonumber\\
&=&  \alpha _{1} \, C  + 
\frac{\alpha _{2}}{n-2}\, E + \frac{\alpha _{8}}{2} \, g \wedge g ,
\end{eqnarray*}
i.e.,
\begin{eqnarray}
\mathrm{Weyl} (T) &=& \alpha _{1} \, C  + \frac{\alpha _{2}}{n-2}\, E 
+ \frac{\alpha _{8}}{2} \, g \wedge g ,
\label{2023.03.06e}
\end{eqnarray}
where $\alpha_{7}$ and $\alpha_{8}$ are some real numbers. 
From (\ref{2023.03.06e}),
by suitable contraction, we get im\-me\-dia\-te\-ly 
$\alpha _{8} = 0$, and in a consequence (\ref{2023.03.06c}),  
completing the proof.

\section{Pseudosymmetry type curvature conditions}

It is well-known that if a semi-Riemannian manifold $(M,g)$, 
$\dim M = n \geq 3$, 
is locally symmetric then 
$\nabla R = 0$
on $M$
(see, e.g., {\cite[Chapter 1.5] {Lumiste}}). 
This  
implies  the following integrability condition
${\mathcal{R}}(X,Y ) \cdot R = 0$ 
in short $R \cdot R = 0$.
Semi-Riemannian manifold satisfying the last condition 
is called {\sl semisymmetric} (see, e.g., 
{\cite[Chapter 8.5.3] {TEC_PJR_2015}}, {\cite[Chapter 20.7] {Chen-2011}}, 
{\cite[Chapter 1.6] {Lumiste}}, \cite{{Sz 1}, 
{LV3-Foreword}}).
Semisymmetric manifolds form a subclass of the class 
of pseudosymmetric manifolds.
A semi-Riemannian manifold $(M,g)$, $\dim M = n \geq 3$,
 is said to be {\sl pseudosymmetric} 
if the tensors $R \cdot R$ and $Q(g,R)$ 
are linearly dependent at every point of $M$
(see, e.g., 
{\cite[Chapter 8.5.3] {TEC_PJR_2015}}, 
{\cite[Chapter 20.7] {Chen-2011}},
{\cite[Section 15.1] {Chen-2021}},
{\cite[Chapter 6] {DHV2008}},
{\cite[Chapter 12.4] {Lumiste}}, 
\cite{{DGHHY}, {DGHSaw}, {DHV2008}, {DVV1991}, 
{HV_2007}, {HaVerSigma}, 
{SDHJK}, {LV1}, {LV2}, 
{LV3-Foreword}, {LV4}} and references therein). 
This is equivalent to
\begin{eqnarray}
R \cdot R &=& L_{R}\, Q(g,R) 
\label{pseudo}
\end{eqnarray}
on   ${\mathcal{U}}_{R}  \subset M$,
where $L_{R}$ is some function on ${\mathcal{U}}_{R}$. 
Every semisymmetric manifold is pseudosymmetric.
The converse statement is not true (see, e.g., \cite{{DVV1991}}).
We note that (\ref{pseudo}) implies
\begin{eqnarray}
R \cdot S &=& L_{R}\, Q(g,S) 
\label{Weyl-pseudo-bis}
\end{eqnarray}
and
\begin{eqnarray}
R \cdot C &=&  L_{R}\, Q(g,C) .  
\label{Weyl-pseudo-ter}
\end{eqnarray}
Conditions 
(\ref{pseudo}), (\ref{Weyl-pseudo-bis}) and (\ref{Weyl-pseudo-ter})
are equivalent on the set ${\mathcal{U}}_{S} \cap {\mathcal{U}}_{C}$ 
of any warped product manifold  $\overline{M} \times _{F} \widetilde{N}$, 
with $\dim \overline{M} = \dim \widetilde{N} = 2$ 
(see, e.g., \cite{{2016_DGJZ}} and references therein).

A semi-Riemannian manifold $(M,g)$, $\dim M = n \geq 3$, 
is called {\sl Ricci-pseudosymmetric} 
if the tensors $R \cdot S$ and $Q(g,S)$ 
are linearly dependent at every point of $M$
(see, e.g., {\cite[Chapter 8.5.3] {TEC_PJR_2015}}, 
{\cite[Section 15.1] {Chen-2021}}, \cite{DGHSaw}).
This is equivalent on ${\mathcal{U}}_{S}$ to 
\begin{eqnarray}
R \cdot S &=& L_{S}\, Q(g,S) , 
\label{Riccipseudo07}
\end{eqnarray}
where $L_{S}$ is some function on ${\mathcal{U}}_{S}$. 
Every warped product manifold $\overline{M} \times _{F} \widetilde{N}$
with a $1$-dimensional manifold  $(\overline{M}, \overline{g})$ and
an $(n-1)$-dimensional Einstein semi-Riemannian manifold 
$(\widetilde{N}, \widetilde{g})$, $n \geq 3$, 
and a warping function $F$, 
is a Ricci-pseudosymmetric manifold,
see, e.g., 
{\cite[Section 1] {Ch-DDGP}}
and
{\cite[Example 4.1] {2016_DGJZ}}.

A semi-Riemannian manifold $(M,g)$, $\dim M = n \geq 4$, 
is said to be {\sl Weyl-pseudo\-sym\-met\-ric} 
if the tensors $R \cdot C$ and $Q(g,C)$ are linearly dependent 
at every point of $M$
\cite{{DGHHY}, {DGHSaw}}. 
This is equivalent on ${\mathcal{U}}_{C}$ to 
\begin{eqnarray}
R \cdot C &=& L_{1}\, Q(g,C) ,  
\label{Weyl-pseudo}
\end{eqnarray}
where $L_{1}$ is some function on ${\mathcal{U}}_{C}$. 
We can easily check that 
on every Einstein manifold $(M,g)$, $\dim M \geq 4$,
(\ref{Weyl-pseudo}) turns into
\begin{eqnarray*}
R \cdot R &=& L_{1}\, Q(g,R) .
\end{eqnarray*}
For a presentation of results on the problem 
of the equivalence of pseudosymmetry, 
Ricci-pseudo\-sym\-met\-ry and Weyl-pseudosymmetry
we refer to {\cite[Section 4] {DGHSaw}}.

A semi-Riemannian manifold $(M,g)$, $\dim M = n \geq 4$, is said to be 
a {\sl manifold with pseudosymmetric Weyl tensor}
({\sl to have a pseudosymmetric conformal Weyl tensor})
if the tensors $C \cdot C$ and $Q(g,C)$ 
are linearly dependent at every point of $M$ 
(see, e.g., {\cite[Section 15.1] {Chen-2021}}, 
\cite{{DGHHY}, {DGHSaw}, {2016_DGJZ}}).
This is equivalent on ${\mathcal U}_{C}$ to 
\begin{eqnarray}
C \cdot C &=& L_{C}\, Q(g,C) ,  
\label{4.3.012}
\end{eqnarray}
where $L_{C}$ is some function on ${\mathcal{U}}_{C}$. 
Every warped product manifold 
 $\overline{M} \times _{F} \widetilde{N}$, 
with $\dim \overline{M} = \dim \widetilde{N} = 2$, 
satisfies (\ref{4.3.012})
(see, e.g., \cite{{DGHHY}, {DGHSaw}, {2016_DGJZ}} and references therein).
Thus in particular,
the Schwarzschild spacetime, the Kottler spacetime
and the Reissner-Nordstr\"{o}m spacetime satisfy (\ref{4.3.012}).
Semi-Riemannian manifolds with pseudosymmetric Weyl tensor 
were studied among others in \cite{{DGHHY}, {DeHoJJKunSh}, {DVY-1994}}.

Warped product manifolds $\overline{M} \times _{F} \widetilde{N}$, 
of dimension $\geq 4$,
satisfying on 
${\mathcal U}_{C} \subset \overline{M} \times _{F} \widetilde{N}$,
the condition 
\begin{eqnarray}
R \cdot R - Q(S,R) &=& L\, Q(g,C) ,  
\label{genpseudo01}
\end{eqnarray}
where $L$ is some function on ${\mathcal{U}}_{C}$,
were studied among others in \cite{49}. In that paper
necessary and sufficient conditions for  
$\overline{M} \times _{F} \widetilde{N}$ 
to be a manifold satisfying (\ref{genpseudo01}) are given.
Moreover, in that paper it was proved that 
any $4$-dimensional warped product manifold 
$\overline{M} \times _{F} \widetilde{N}$, 
with a $1$-dimensional base $(\overline{M},\overline{g})$, 
satisfies (\ref{genpseudo01}) {\cite[Theorem 4.1] {49}}.

We refer to
\cite{{Ch-DDGP}, {D-1991}, {DGHHY}, {DGHSaw}, {DGHZ01}, {2016_DGJZ}, 
{DHV2008}, {DeHoJJKunSh}, {SDHJK}} 
for details on semi-Riemannian manifolds satisfying 
(\ref{pseudo}) and (\ref{Riccipseudo07})-(\ref{genpseudo01}), 
as well other conditions of this kind, named pseudosymmetry 
type curvature conditions. 
We also refer to {\cite[Section 3] {DeHoJJKunSh}} 
for a recent survey on manifolds 
satisfying such curvature conditions.
It seems that the condition (\ref{pseudo}) 
is the most important condition of that family of curvature conditions
(see, e.g., \cite{2016_DGJZ}).
The Schwarzschild spacetime, the Kottler spacetime, 
the Reissner-Nordstr\"{o}m spacetime, 
as well as the Friedmann-Lema{\^{\i}}tre-Robertson-Walker spacetimes 
are the ``oldest'' examples 
of pseudosymmetric warped product manifolds 
(see, e.g., \cite{{2016_DGJZ}, {DHV2008}, {DVV1991}, {SDHJK}}).
We finish this section with the following remarks.

\begin{rem} (i) In view of {\cite[Lemma 3.2 (ii)] {D-1990}},
we can state that the following identity is satisfied 
on every semi-Riemannian manifold 
$(M,g)$, $\dim M = n \geq 3$,
with vanishing Weyl conformal curvature tensor $C$
\begin{eqnarray} 
R \cdot R - Q(S,R) &=&
\frac{1}{(n-2)^{2}}\,
Q(g, g \wedge S^{2} + \frac{n-2}{2} \, S \wedge S - \kappa \, g \wedge S) .
\label{2022.11.20.zz}
\end{eqnarray}
From (\ref{2022.11.20.zz}), by (\ref{2022.11.10.bbb}), we get 
\begin{eqnarray*} 
R \cdot R - Q(S,R) &=& \frac{1}{(n-2)^{2}}\, Q(g,E) ,
\end{eqnarray*}
where the tensor $E$ is defined by  (\ref{2022.11.10.aaa}). In particular, if $n = 3$ then $E = 0$ on $M$.
\newline
(ii) As it was stated in {\cite[Theorem 3.1] {D-1990}} 
on every $3$-dimensional semi-Riemannian manifold $(M,g)$
the identity $R \cdot R = Q(S,R)$ is satisfied.
\newline
(iii) From (i) it follows  that
on every semi-Riemannian conformally flat manifold 
$(M,g)$, $\dim M = n \geq 4$, the conditions: $R \cdot R = Q(S,R)$
and (\ref{2022.11.20.aa}) are equivalent.
\end{rem}

\begin{rem}
Let $(M,g)$, $\dim M = n \geq 4$, be a semi-Riemannian manifold.
\newline
(i) {\cite[Theorem 3.4 (i)] {2016_DGJZ}}
The following identity is satisfied on ${\mathcal{U}}_{C} \subset M$
\begin{eqnarray} 
C \cdot R + R \cdot C &=& R \cdot R + C \cdot C - \frac{1}{ (n-2)^{2}}
Q( g, g \wedge S^{2} - \frac{\kappa }{n-1}\, g \wedge S) .
\label{2022.12.29.aaa}
\end{eqnarray}
(ii) If 
(\ref{genpseudo01}) holds on ${\mathcal{U}}_{C} \subset M$ 
then (\ref{2022.12.29.aaa}) turns into 
\begin{eqnarray}
& & 
C \cdot R + R \cdot C \ =\ C \cdot C + Q(S,R) + L  \, Q(g,C)\nonumber\\
& & 
- \frac{1}{ (n-2)^{2}}\,
Q( g, g \wedge S^{2} + \frac{n-2}{2}\, S \wedge S 
- \frac{\kappa }{n-1}\, g \wedge S) .
\label{2022.12.29.ddd}
\end{eqnarray}
Moreover, from (\ref{2022.12.29.ddd}), by an application 
of (\ref{abcd07}) and (\ref{2022.11.10.bbb}), 
we get on ${\mathcal{U}}_{C} \subset M$
\begin{eqnarray}
C \cdot R + R \cdot C
&=& C \cdot C + Q(S,C) + L \, Q(g,C) - \frac{1}{ (n-2)^{2}} \, Q( g, E) ,
\label{2022.12.29.fff}
\end{eqnarray}
where the tensor $E$ is defined by  (\ref{2022.11.10.aaa}).
\newline
(iii) (cf. {\cite[Theorem 3.4 (iii)] {2016_DGJZ}})
If (\ref{4.3.012}) and (\ref{genpseudo01}) 
hold on ${\mathcal{U}}_{C} \subset M$ then
(\ref{2022.12.29.fff}) turns into 
\begin{eqnarray*}
C \cdot R + R \cdot C
&=& Q(S,C) + (L_{C} + L) \, Q(g,C) - \frac{1}{ (n-2)^{2}} \, Q( g, E) .
\end{eqnarray*}
\end{rem}

\section{Roter spaces}

Some results of \cite{ {DGHHY}, {DH-1998}, {DY-1994}} 
(cf. {\cite[Section 1] {2016_DGJZ}})
we can present in the following proposition.
\begin{prop}
Let $(M,g)$, $\dim M = n \geq 4$,
be a non-conformally flat and non-Einstein semi-Riemannian manifold.
\newline 
(i) {\cite[Theorem 3.1, Teorem 3.2 (ii)] {DY-1994}}
If (\ref{pseudo}) and (\ref{4.3.012})
hold on ${\mathcal{U}}_{S} \cap {\mathcal{U}}_{C} \subset M$ then
at every point $x \in {\mathcal{U}}_{S} \cap {\mathcal{U}}_{C}$
(\ref{quasi02}) or (\ref{eq:h7a}) is satisfied.
\newline
(ii) {\cite[Theorem 3.1, Teorem 3.2 (ii)] {DH-1998}}
If (\ref{pseudo}) and (\ref{genpseudo01})
hold on ${\mathcal{U}}_{S} \cap {\mathcal{U}}_{C} \subset M$ then
at every point $x \in {\mathcal{U}}_{S} \cap {\mathcal{U}}_{C}$
(\ref{quasi02}) or (\ref{eq:h7a}) is satisfied.
\newline
(iii) ({cf. \cite[Proposition 3.2, Theorem 3.3, Theorem 4.4] {DGHHY}}) 
If (\ref{4.3.012}), (\ref{genpseudo01}) 
and $R \cdot S = Q(g,D)$, 
for some symmetric $(0,2)$-tensor $D$,
hold on ${\mathcal{U}}_{S} \cap {\mathcal{U}}_{C} \subset M$ then 
at every point $x \in {\mathcal{U}}_{S} \cap {\mathcal{U}}_{C}$
(\ref{quasi02}) or (\ref{eq:h7a}) is satisfied.
\end{prop}

We recall that
a  non-quasi-Einstein and non-conformally flat 
semi-Riemannian manifold $(M,g)$, $\dim M = n \geq 4$, 
satisfying (\ref{eq:h7a}) on 
${\mathcal U}_{S} \cap {\mathcal U}_{C} \subset M$ 
is called a {\sl Roter type manifold}, or a {\sl Roter manifold}, 
or a {\sl Roter space} 
(see, e.g., {\cite[Section 15.5] {Chen-2021}},
\cite{{P106}, {2016_DGJZ}, {DGP-TV02}, {DHV2008}}).

Roter spaces and in particular Roter hypersurfaces 
in semi-Riemannian spaces of constant curvature were studied in:
\cite{{DDH-2021}, {P106}, {DGHHY}, {DGHZ01},
{R102}, {DeKow}, {DePlaScher}, {DeScher}, {G5}, {Kow01}, {Kow02}}. 
In particular, 
(\ref{pseudo}) and (\ref{Riccipseudo07})-(\ref{genpseudo01})
are satisfied on such manifolds. More precisely, we have

\begin{thm} (see, e.g., \cite{{DGHSaw}, {2016_DGJZ}}, 
{\cite[eq. (28)]{DP-TVZ}}) 
If $(M,g)$, $\dim M = n \geq 4$, is a semi-Riemannian 
Roter space satisfying
(\ref{eq:h7a}) on ${\mathcal U}_{S} \cap {\mathcal U}_{C} \subset M$ 
then on this set we have:
\begin{eqnarray*}
S^{2} &=& \alpha _{1}\, S + \alpha _{2} \, g ,
\ \ \
\alpha _{1} 
\ =\ 
\kappa + \frac{(n-2)\mu -1 }{\phi} ,
\ \ \
\alpha _{2}
\ =\
\frac{\mu \kappa + (n-1) \eta }{\phi } ,\\ 
R \cdot C &=& L_{R}\, Q(g,C),
\ \ \
L_{R} \ =\ \frac{1}{\phi} 
\left(  (n-2) (\mu ^{2} - \phi \eta) - \mu \right) ,\\
R \cdot R & = & L_{R}\, Q(g,R) ,\\
R \cdot S & = & L_{R}\, Q(g,S),
\end{eqnarray*}
\begin{eqnarray*}
R \cdot R &=& Q(S,R) + L \, Q(g,C) ,
\ \ \
L \ =\ L_{R} + \frac{\mu }{\phi }
\ =\
\frac{n-2}{\phi} (\mu ^{2} - \phi \eta),\\
C \cdot C &=& L_{C}\, Q(g,C) ,
\ \ \ 
L_{C} \ =\ L_{R} + \frac{1}{n-2} (\frac{\kappa }{n-1} - \alpha _{1} ) ,\\
C \cdot R & = & L_{C}\, Q(g,R),\\
C \cdot S & = & L_{C}\, Q(g,S) ,
\end{eqnarray*}
\begin{eqnarray*}  
C \cdot R + R \cdot C 
&=& Q(S,C) + \left( L + L_{C} - \frac{1}{ (n-2) \phi } \right) Q(g,C) ,\\
R \cdot C - C \cdot R &=& 
\left( \frac{1}{\phi} ( \mu - \frac{1}{n-2} ) 
+ \frac{\kappa }{n-1} \right) Q(g,R)
+ 
\left( \frac{\mu}{\phi } ( \mu - \frac{1}{n-2}) - \eta \right) Q(S,G) ,
\end{eqnarray*}
(\ref{Roterformula}) and 
\begin{eqnarray*}  
R\cdot C - C \cdot R 
&=& 
Q \left( \left(\frac{\mu \kappa  }{n-1} + \eta \right) \, g 
+ 
\left( \frac{1 }{n-2} - \mu - \frac{\phi \kappa }{n-1} \right)\, S, 
g \wedge S \right) . 
\end{eqnarray*}
\end{thm}

\begin{rem}
(i) In the standard Schwarzschild coordinates 
$(t; r; \theta; \phi)$, 
and the physical units ($c = G = 1$), 
the Reissner-Nordstr\"{o}m-de Sitter ($\Lambda > 0$), and 
Reissner-Nordstr\"{o}m-anti-de Sitter ($\Lambda < 0$) 
spacetimes are given by the line element (see, e.g., \cite{SH})
\begin{eqnarray}
ds^{2} &=& - h(r)\, dt^{2} + h(r)^{-1}\, dr^{2} + r^{2}\, 
( d\theta^{2} + \sin ^{2}\theta \, d\phi^{2}) ,\label{rns01}\\
h(r) &=& 1 - \frac{2 M}{r} + \frac{Q^{2}}{r^{2}} 
- \frac{\Lambda }{3} r^{3} ,\nonumber
\end{eqnarray}
where 
$M$, $Q$ and $\Lambda$ are non-zero constants. 
\newline
(ii) 
{\cite[Section 6] {2017_DGHZ}} 
(see also
{\cite[Remark 2 (ii)] {DDH-2021}},
{\cite[Remark 2.1 (ii)] {DH-Colloq}})
The metric (\ref{rns01}) satisfies (\ref{eq:h7a}) with 
\begin{eqnarray*}
\phi &=&
\frac{3}{2} ( Q^{2}  -  M r) r^{4} }{2 Q^{-4},\ \ \ 
\mu \ =\
\frac{1}{2} (Q^{4} + 3 Q^{2} \Lambda r^{4} - 3 \Lambda M r^{5}) Q^{-4}, \\ 
\eta
&=&
\frac{1}{12}
( 3 Q^{6}
+ 4 Q^{4} \Lambda r^{4}
- 3 Q^{4} M r
+ 9 Q^{2} \Lambda ^{2} r^{8}
- 9 \Lambda^{2} M r^{9})  r^{-4} Q^{-4} .
\end{eqnarray*}
If we set $\Lambda = 0$ in (\ref{rns01}) then 
we obtain the line element of the Reissner-Nordstr\"{o}m
spacetime, see, e.g., {\cite[Section 9.2] {GrifPod}} 
and references therein. 
It seems that the Reissner-Nordstr\"{o}m spacetime is the oldest example 
of the Roter warped product space.
\newline
(iii) 
In \cite{DH-Colloq}
a particular class of Roter warped product spaces was determined 
such that every manifold of that class 
admits a non-trivial geodesic mapping onto some Roter warped product space.
Moreover, both geodesically related manifolds 
are pseudosymmetric of constant type. 
\newline
(iv)
An algebraic classification of the Roter type $4$-dimensional 
spacetimes is given in
\cite{DDH-2021}.
\newline
(v)
Some comments on pseudosymmetric manifolds (also called 
Deszcz symmetric spaces),
as well as Roter spaces, are given in {\cite[Section 1] {DecuP-TSVer}}
(see also {\cite[Remark 2 (iii)] {DDH-2021}},
{\cite[Remark 2.1 (iii)] {DH-Colloq}}):
"{\sl{From a geometric point of view, the Deszcz symmetric spaces may well
be considered to be the simplest Riemannian manifolds next 
to the real space forms.}}" 
and 
"{\sl{From an algebraic point of view, Roter spaces may well 
be considered to
be the simplest Riemannian manifolds next to the real space forms.}}"
For further comments we refer to \cite{LV3-Foreword}.
\end{rem}

We finish this section with the following results. 
\begin{prop} {\cite[Lemma 2.2] {DGHSaw-2022}}
If $(M,g)$, $\dim M = n \geq 4$, is a Roter space satisfying
(\ref{eq:h7a}) 
on ${\mathcal U}_{S} \cap  {\mathcal U}_{C} \subset M$ 
then (\ref{2022.11.22.aa}) holds 
on this set, i.e., $(n-2)\, C = \phi \, E$, 
where the tensor $E$ is defined by  (\ref{2022.11.10.aaa}).
\end{prop}

Propositions 2.1, 4.1 and 4.4 lead to the following
\begin{prop}
Let $(M,g)$, $\dim M = n \geq 4$,
be a non-conformally flat and non-Einstein semi-Riemannian manifold.
If 
(\ref{pseudo}) and (\ref{4.3.012}),
or
(\ref{pseudo}) and (\ref{genpseudo01}),
or
(\ref{4.3.012}), (\ref{genpseudo01}) 
and $R \cdot S = Q(g,D)$, 
for some symmetric $(0,2)$-tensor $D$,
hold on ${\mathcal{U}}_{S} \cap {\mathcal{U}}_{C} \subset M$ then
$E = \lambda \, C$
on ${\mathcal{U}}_{S} \cap {\mathcal{U}}_{C}$,
where the tensor $E$ is defined by  (\ref{2022.11.10.aaa}) 
and $\lambda $ is some function on this set. 
\end{prop}

\section{Warped product manifolds with 2-dimensional base manifold}

Proposition 2.1 and Proposition 4.4 imply  
\begin{prop} {\cite[Proposition 2.3] {DGHSaw-2022}}
If $(M,g)$, $\dim M = n \geq 4$, is a semi-Riemannian manifold satisfying
(\ref{quasi02}) or (\ref{eq:h7a}) at every point 
of ${\mathcal U}_{S} \cap  {\mathcal U}_{C} \subset M$ 
then the following equation is satisfied on this set  
\begin{eqnarray}
\tau \, C &=&  
g \wedge S^{2} + \frac{n-2}{2} \, S \wedge S - \kappa \, g \wedge S
+ \frac{\kappa ^{2} - \mathrm{tr}_{g} (S^{2})}{2(n-1)}
\, g \wedge g ,
\label{2022.11.22.jj}
\end{eqnarray}
where $\tau$ is some function 
on ${\mathcal U}_{S} \cap {\mathcal U}_{C}$.
\end{prop}

Proposition 5.1, {\cite[Theorem 7.1 (ii)] {2016_DGJZ}} and
{\cite[Theorem 4.1] {DeKow}} imply 
\begin{thm} {\cite[Theorem 2.4] {DGHSaw-2022}}
Let $\overline{M} \times _{F} \widetilde{N}$ 
be the warped product manifold 
with a $2$-dimensional semi-Riemannian manifold 
$(\overline{M},\overline{g})$,
an $(n-2)$-dimensional semi-Riemannian manifold 
$(\widetilde{N},\widetilde{g})$, 
$n \geq 4$, a warping function $F$, 
and let $(\widetilde{N},\widetilde{g})$ 
be a space of constant curvature when $n \geq 5$. 
Then (\ref{2022.11.22.jj}) holds on
${\mathcal U}_{S} \cap {\mathcal U}_{C} 
\subset \overline{M} \times _{F} \widetilde{N}$.
\end{thm}

\begin{example} {\cite[Example 2.1] {DGHSaw-2022}}
(i) Let $\mathbb{S}^{k}(r)$ be a $k$-dimensional standard sphere 
of radius $r$ in $\mathbb{E}^{k+1}$, $k \geq 1$.
It is well-known that 
the Cartesian product 
$\mathbb{S}^{1}(r_{1}) \times \mathbb{S}^{n-1}(r_{2})$
of spheres $\mathbb{S}^{1}(r_{1})$ and $\mathbb{S}^{n-1}(r_{2})$, 
$n \geq 4$,
and more generally,
the warped product manifold 
$\mathbb{S}^{1}(r_{1}) \times _{F} \mathbb{S}^{n-1}(r_{2})$
of spheres $\mathbb{S}^{1}(r_{1})$ and $\mathbb{S}^{n-1}(r_{2})$,
$n \geq 4$, with a warping function $F$,
is a conformally flat manifold.
\newline
(ii)
As it was stated in {\cite[Example 3.2] {G5}},
the Cartesian product 
$\mathbb{S}^{p}(r_{1}) \times \mathbb{S}^{n-p}(r_{2})$
of spheres 
$\mathbb{S}^{p}(r_{1})$ and $\mathbb{S}^{k}(r_{2})$ 
such that 
$2 \leq p \leq n-2$ and
$(n-p-1) r_{1}^{2} \neq (p-1) r_{2}^{2}$
is a non-conformally flat and non-Einstein manifold
satisfying the Roter equation (\ref{eq:h7a}) on 
${\mathcal U}_{S} \cap  {\mathcal U}_{C} 
= \mathbb{S}^{p}(r_{1}) \times \mathbb{S}^{n-p}(r_{2})$. 
\newline
(iii) 
{\cite[Example 4.1] {DeKow}}
The warped product manifold 
$\mathbb{S}^{p}(r_{1}) \times _{F} \mathbb{S}^{n-p}(r_{2})$,
$2 \leq p \leq n-2$, with some special warping function $F$,
satisfies on ${\mathcal U}_{S} \cap  {\mathcal U}_{C} 
\subset \mathbb{S}^{p}(r_{1}) \times \mathbb{S}^{n-p}(r_{2})$ 
the Roter equation (\ref{eq:h7a}). Thus some 
warped product manifolds 
$\mathbb{S}^{2}(r_{1}) \times _{F} \mathbb{S}^{n-2}(r_{2})$ 
are Roter spaces.
\newline
(iv)  Properties of pseudosymmetry type of warped products 
with $2$-dimensional base manifold,
a warping function $F$,
and an $(n-2)$-dimensional fibre, 
$n \geq 4$,
assumed to be of constant curvature when $n \geq 5$, 
were determined in  
{\cite[Sections 6 and 7] {2016_DGJZ}}.
Evidently, 
warped product manifolds 
$\mathbb{S}^{2}(r_{1}) \times _{F} \mathbb{S}^{n-2}(r_{2})$,
$n \geq 4$, are such manifolds.
Let $g$, $R$, $S$, $\kappa$ and $C$ denote the metric tensor, 
the Riemann-Christoffel curvature tensor, the Ricci tensor,
the scalar curvature and  
the Weyl conformal curvature tensor of 
$\mathbb{S}^{2}(r_{1}) \times _{F} \mathbb{S}^{n-2}(r_{2})$, respectively.
From {\cite[Theorem 7.1] {2016_DGJZ}}
it follows that on set $V$ of all points of 
${\mathcal U}_{S} \cap  {\mathcal U}_{C} 
\subset \mathbb{S}^{2}(r_{1}) \times _{F} \mathbb{S}^{n-2}(r_{2})$
at which the tensor $S^{2}$ is not a linear combination 
of the tensors $S$ and $g$,
the Weyl tensor $C$ is expressed by
\begin{eqnarray}
C &=& \lambda \left(
g \wedge S^{2} + \frac{n-2}{2} \, S \wedge S - \kappa \, g \wedge S
+ \frac{\kappa ^{2} - \mathrm{tr}_{g} (S^{2}) }{2(n-1)} 
\, g \wedge g \right) ,
\label{2022.11.10.zz}
\end{eqnarray}
where $\lambda$ is some function on $V$.
This,  by (\ref{Weyl}), turns into
\begin{eqnarray*}
R &=& \lambda \,
g \wedge S^{2} + \frac{n-2}{2} \lambda \, S \wedge S 
+ \left( \frac{1}{n-2} - \kappa \lambda \right) g \wedge S\\
& &
+ \frac{1}{2(n-1)} \left( (\kappa ^{2} - \mathrm{tr}_{g} (S^{2}) ) \lambda
- \frac{\kappa }{n-2} \right) g \wedge g .
\end{eqnarray*}
Thus  (\ref{B001simply}) is satisfied on $V$.
Moreover,  (\ref{eq:h7a}) holds 
at all points of $({\mathcal U}_{S} \cap  {\mathcal U}_{C}) \setminus V$,
at which (\ref{quasi02}) is not satisfied.
From Lemma 2.2 it follows that (\ref{2022.11.10.zz})
holds at all points of ${\mathcal U}_{S} \cap  {\mathcal U}_{C} \subset
\mathbb{S}^{2}(r_{1}) \times _{F} \mathbb{S}^{n-2}(r_{2})$, $n \geq 4$,
at which (\ref{quasi02}) is not satisfied. Finally, in view of Theorem 5.2, 
we can state that (\ref{2022.11.22.jj}) holds on
${\mathcal U}_{S} \cap  {\mathcal U}_{C}$.
\end{example}

\section{Essentially conformally symmetric manifolds}

Let $(M,g)$, $\dim M = n \geq 4$, be a semi-Riemannian manifold 
with parallel Weyl conformal curvature tensor,
i.e. $\nabla C  = 0$ on $M$. 
It is obvious that the last condition implies $R \cdot C  = 0$.
Moreover, let the manifold $(M,g)$ be neither conformally flat nor 
locally symmetric. 
Such manifolds are called {\sl{essentially conformally 
symmetric manifolds}}, e.c.s. manifolds/metrics, or ECS manifolds/metrics, 
in short (see, e.g.,
\cite{{DerRot01}, {DerRot02}, {DerRot2007PADGE}, 
{DerTer2210a}, {DerTer2210b}, {Hotlos}}).
E.c.s. manifolds 
are semisymmetric manifolds
($R \cdot R  = 0$, {\cite[Theorem 9] {DerRot01}})
satisfying
$\kappa  = 0$ and $Q(S,C) \, =\, 0$
({\cite[Theorems 7 and 8] {DerRot01}}).
In addition, 
\begin{eqnarray}
F\, C &=& \frac{1}{2}\, S \wedge S 
\label{roter77}
\end{eqnarray}
holds on $M$, where $F$ is some function on $M$, called the 
{\sl{fundamental function}} \cite{DerRot02}.  
At every point of $M$ we also have $\mbox{rank}\, S \leq 2$ 
{\cite[Theorem 5] {DerRot02}}.
We mention that the local structure 
of e.c.s. manifolds is already determined.
We refer to \cite{{DerRot2007}, {DerRot2009}} for the final results 
related to this subject.
We also mention that certain e.c.s. metrics are realized 
on compact manifolds 
(\cite{{DerRot2008}, {DerRot2010}, {DerTer2210a}, {DerTer2210b}}).

Equation (\ref{roter77}), by suitable contraction, leads immediately to
$S^{2} = \kappa \, S$, which by $\kappa = 0$, reduces to $S^{2} = 0$. 
Evidently, $\mathrm{tr}_{g} (S^{2}) = 0$.
Now using (\ref{roter77}) we get (\ref{2023.03.03.a}).
Thus we have
\begin{thm}
Condition 
(\ref{2022.11.22.jj}), with $\tau = (n-2) F$,
is satisfied on every essentially conformally symmetric manifold $(M,g)$.
\end{thm}

\section{Hypersurfaces in semi-Riemannian conformally flat spaces}

Let $M$, $\dim M = n \geq 4$, be a hypersurface 
isometrically immersed in a semi-Riemannian conformally flat manifold 
$N$, $\dim N = n + 1$.
Let $g_{ad}$, $H_{ad}$, 
$G_{abcd} = g_{ad}g_{bc} - g_{ac}g_{bd}$ and $C_{abcd}$
be the local components of the metric tensor $g$, 
the second fundamental tensor $H$, the $(0,4)$-tensor $G$ 
and the Weyl conformal curvature tensor $C$ of $M$, 
respectively. 
As it was stated in {\cite[eq. (20)] {DV-1991}} 
(see also {\cite[eq. (11)] {P38}}) we have 
\begin{eqnarray}
C_{abcd} 
&=& 
\varepsilon \, (H_{ad}H_{bc} - H_{ac}H_{bd}) 
- 
\frac{\varepsilon \, \mathrm{tr} (H)}{ n-2}\, 
( g_{ad}H_{bc} + g_{bc}H_{bd} - g_{ac}H_{bd} - g_{bd}H_{ac} ) 
\nonumber\\
& &
+ 
\frac{\varepsilon }{ n-2}\, 
( g_{ad}H^{2}_{bc} + g_{bc}H^{2}_{bd} - g_{ac}H^{2}_{bd} - g_{bd}H^{2}_{ac} ) 
+ 
\mu \, G_{abcd} \, ,  
\label{2022.08.03.bb}
\end{eqnarray}
where
$\varepsilon = \pm 1$, 
$\mathrm{tr} (H) = g^{ad}H_{ad}$,  
$H^{2}_{ad} = g^{bc}H_{ab}H_{cd}$
and $\mu$ is some function on $M$.
From (\ref{2022.08.03.bb}), by contraction we get easily
\begin{eqnarray}
\mu 
&=& \frac{ \varepsilon }{ (n-2)(n-1)}\, 
( (\mathrm{tr} (H))^{2} - \mathrm{tr} (H^{2}) ) \, ,
\label{2022.08.03.cc}
\end{eqnarray}
where  
$\mathrm{tr} (H^{2}) 
\, =\, g^{ad}H^{2}_{ad}$. 
Now (\ref{2022.08.03.bb}) and (\ref{2022.08.03.cc}) yield
\begin{eqnarray}
C &=& \frac{\varepsilon}{n-2} 
\left( g \wedge H^{2} + \frac{n-2}{2}\, H \wedge H 
- \mathrm{tr} (H) \, g \wedge H 
+ \frac{(\mathrm{tr} (H))^{2} 
- \mathrm{tr}(H^{2})}{2 (n-1)} \, g \wedge g \right).
\label{2022.08.03.dd}
\end{eqnarray}
If $H = \frac{\mathrm{tr} (H)}{n} g$ at a point $x \in M$,
i.e., $M$ is umbilical at $x$, then from (\ref{2022.08.03.dd}) it follows
immediately that the tensor $C$ vanishes at $x$.
If at a non-umbilical point $x \in M$, we have 
$\mathrm{rank} (H - \alpha g) = 1$,
for some $\alpha \in \mathbb{R}$, 
i.e., $M$ is quasi-umbilical at $x$,
then in view of Proposition 2.1 (i), the tensor $C$ vanishes at $x$. 
Conversely, if at a non-umbilical point $x \in M$ the tensor $C$ vanishes 
then in view of Proposition 2.1 (ii) we have $\mathrm{rank} (H - \alpha g) = 1$,
for some $\alpha \in \mathbb{R}$. 
Thus we can present {\cite[Theorem 4.1] {DV-1991}} in the folowing form.
\begin{thm} 
Let $M$, $\dim M = n \geq 4$, be a hypersurface 
isometrically immersed in a semi-Riemannian conformally flat manifold 
$N$, $\dim N = n + 1$.
At every non-umbilical point $x \in M$ the tensor $C$ vanishes at $x$
 if and only if
at $x$ we have $\mathrm{rank} (H - \alpha g) = 1$,
for some $\alpha \in \mathbb{R}$. 
\end{thm}

\begin{rem}
Let $M$, $\dim M = n \geq 4$, be a hypersurface 
isometrically immersed in a semi-Riemannian conformally flat manifold 
$N$, $\dim N = n + 1$.
\newline
(i)
We assume that at all points of ${\mathcal U}_{C} \subset M$ the tensor $H^{2}$ 
is a linear combination of $H$ and $g$, i.e., 
\begin{eqnarray}
H^{2} &=& \alpha_{1}\, H + \alpha_{2}\, g 
\label{2022.12.12.mm}
\end{eqnarray}
on ${\mathcal U}_{C}$,
where $\alpha_{1}$ and $\alpha_{2}$ are some functions on this set. 
Now (\ref{2022.08.03.dd}) turns into
\begin{eqnarray}
C 
&=& 
\frac{\varepsilon }{2}\, H \wedge H 
+ \frac{\varepsilon (\alpha_{1} - \mathrm{tr} (H) ) }{n-2} 
 \, g \wedge H 
+ \frac{\varepsilon}{n-2}
\left( \alpha_{2} +  \frac{(\mathrm{tr} (H))^{2} 
- \mathrm{tr}(H^{2})}{2 (n-1)} \right)  g \wedge g \nonumber\\
&=& 
\frac{\alpha }{2}\, H \wedge H 
+ \beta \, g \wedge H 
+ \frac{\gamma }{2} \, g \wedge g ,
\label{2022.08.03.ff}
\end{eqnarray}
where
\begin{eqnarray}
\alpha \ =\ \varepsilon,\ \ \ 
\beta \ =\  \frac{\varepsilon (\alpha_{1} - \mathrm{tr} (H) ) }{n-2} ,\ \ \
\gamma \ =\ \frac{\varepsilon}{n-2}
\left( 2 \alpha_{2} +  \frac{(\mathrm{tr} (H))^{2} 
- \mathrm{tr}(H^{2})}{n-1} \right) . 
\label{2022.08.03.gg}
\end{eqnarray}
From (\ref{2022.08.03.ff}) and (\ref{2022.08.03.gg}),
in view of {\cite[Theorem 3.1 (i)] {Kow02}}, we get
\begin{eqnarray}
C \cdot C 
&=& (n-2) \left( \frac{\beta ^{2}}{\alpha } - \gamma \right)  Q(g,C) 
\ =\  (n-2) ( \varepsilon \beta ^{2} - \gamma )\,  Q(g,C)
\label{2022.08.03.hh}
\end{eqnarray}
on ${\mathcal U}_{C}$, with $\alpha$, $\beta$ and $\gamma$ 
defined by (\ref{2022.08.03.gg}).
Thus $M$ is a hypersurface with pseudosymmetric Weyl tensor.
We also note that from (\ref{2022.12.12.mm}) we get immediately
$\alpha _{2} = \frac{1}{n} ( \mathrm{tr}(H^{2}) - \alpha _{1} \mathrm{tr} (H) )$
and 
\begin{eqnarray*}
H^{2} - \frac{\mathrm{tr}(H^{2})}{n}\, g 
&=& \alpha_{1} \left( H - \frac{\mathrm{tr} (H)}{n}\, g \right) .
\end{eqnarray*}
(ii) The above presented result, i.e., 
if (\ref{2022.12.12.mm}) is satisfied at every point 
of ${\mathcal U}_{C} \subset M$
then (\ref{4.3.012}) holds on this set, was already obtained in 
{\cite[Proposition 3.1] {DVY-1999}}.
We mention that Proposition 3.1 of \cite{DVY-1999} 
was proved without application of {\cite[Theorem 3.1 (i)] {Kow02}}.
\newline
(iii) 
We assume that the tensor $H$ satisfies on ${\mathcal U}_{C} \subset M$
\begin{eqnarray}
H^{3} &=& \mathrm{tr} (H) \,  H^{2} + \psi \,   H,  
\label{2022.10.01.kk}
\end{eqnarray}
where $\psi $ is some function on this set
and the $(0,2)$-tensor $H^{3}$ is defined by 
$H^{3}_{ad} = g^{bc}H^{2}_{ab}H_{cd}$. 
Then 
\begin{eqnarray}
C \cdot C &=& 
\left(
\frac{\varepsilon }{(n-2)(n-1)} 
(
(\mathrm{tr} (H))^{2} - \mathrm{tr}(H^{2})
)
+ \frac{\varepsilon \psi}{n-2} 
\right)
Q(g,C)\nonumber\\
& &
- \frac{n-3 }{n-2} \, Q(H^{2}, \frac{1}{2}\, H \wedge H) 
\label{2022.10.01.ll}
\end{eqnarray}
on ${\mathcal U}_{C}$ {\cite[eq. (10)] {P37}}, see also 
{\cite[the proof of Lemma 4.1] {DVY-1999}}. 
We refer to \cite{{P37}, {DVY-1999}} for further results on hypersurfaces 
$M$ in conformally flat manifold $N$ satisfying (\ref{2022.10.01.kk}).
\newline
(iv) Recently curvature properties of pseudosymmetry type of hypersurfaces 
isometrically immersed in a semi-Riemannian conformally flat manifold
were investigated in \cite{SE} and \cite{Kundu}.  
\end{rem}

\section{Hypersurfaces in semi-Riemannian space forms}

Let now $N_{s}^{n+1}(c)$, $n  \geq 4$, 
be a semi-Riemannian space of constant curvature 
with signature $(s,n+1-s)$, where
$c = \frac{\widetilde{\kappa}}{n (n+1)}$ and
$\widetilde{\kappa}$ is its scalar curvature. 
Let $M$, $\dim M = n \geq 4$, be a connected hypersurface
isometrically immersed in $N_{s}^{n+1}(c)$.
We denote by
$g$, $R$, $S$, $\kappa$ and $C$, 
the metric tensor, the Riemann-Christoffel curvature tensor, 
the Ricci tensor,
the scalar curvature and the Weyl conformal curvature tensor 
of the hypersurface $M$, respectively.
The Gauss equation of $M$ in $N_{s}^{n+1}(c)$ reads 
(see, e.g., \cite{{2016_DGHZhyper}, {DGP-TV02}, 
{DGPSS}, {2020_DGZ}, {Saw114}})
\begin{eqnarray}
R - \frac{\widetilde{\kappa}}{2 n (n+1)}\, g \wedge g
&=& 
\frac{\varepsilon }{2}  \,  H \wedge H ,\  \  \  \varepsilon \ =\ \pm 1 .
\label{realC5}
\end{eqnarray}  
From (\ref{realC5}), by suitable contractions, we obtain
\begin{eqnarray}
S -   \frac{ (n-1) \widetilde{\kappa}}{ n (n+1)}\, g
&=&  \varepsilon \, ( \mathrm{tr} (H)\, H - H^{2}) ,
\label{realC6}\\
\frac{ \kappa }{n-1}  - \frac{  \widetilde{\kappa} }{n+1}   
&=& \frac{ \varepsilon}{n-1} 
\, ( ( \mathrm{tr} (H))^{2} - \mathrm{tr} (H^{2}) ) .
\label{realC6bb}
\end{eqnarray}
Now using (\ref{realC5}),  (\ref{realC6}) and (\ref{realC6bb}) 
we get immediately
\begin{eqnarray}
Q(H^{2}, \frac{1}{2}\, H \wedge H)
&=&
- Q(  \mathrm{tr} (H)\, H  - H^{2} , \frac{1}{2}\, H \wedge H)  \nonumber\\
&=&
- Q( \varepsilon ( \mathrm{tr} (H)\, H  -  H^{2} ), 
\frac{\varepsilon}{2}\, H \wedge H) \nonumber\\
&=&
- Q\left( S - \frac{ (n-1) \widetilde{\kappa}}{ n (n+1)}\, g ,  
R - \frac{\widetilde{\kappa}}{2 n (n+1)}\, g \wedge g \right ) .
\label{realRa}
\end{eqnarray}

We also recall that the curvatere condition of pseudosymmetry type
(\ref{genpseudo01}) is satisfied on $M$. Precisely, we have on $M$
{\cite[Proposition 3.1] {DV-1991}} 
(see also {\cite[eqs. (3.3) and (3.4)] {2020_DGZ}})
\begin{eqnarray}
R \cdot R - Q(S,R) &=& - \frac{ (n-2) \widetilde{\kappa} }{n (n+1)}\, Q(g,C) .
\label{2023.03.21.aa}
\end{eqnarray}
Now (\ref{2022.12.29.fff}), by (\ref{2023.03.21.aa}), turns into
({\cite[Propopsition 4.7, eq. (4.36)] {2020_DGZ}})
\begin{eqnarray}
C \cdot R + R \cdot C
&=& C \cdot C + Q(S,C) 
- \frac{ (n-2) \widetilde{\kappa} }{n (n+1)}\,  Q(g,C) 
- \frac{1}{ (n-2)^{2}} \, Q( g, E) .
\label{2022.12.29.rrr}
\end{eqnarray}

Let ${\mathcal U}_{H} \subset M$ be the set of all points 
at which the tensor 
$H^{2}$ is not a linear combination of the metric tensor $g$ 
and the second fundamental tensor $H$ of $M$. We have
${\mathcal U}_{H} \subset {\mathcal U}_{S} \cap {\mathcal U}_{C} \subset M$ 
(see, e.g., \cite{{R99}, {DGPSS}, {2020_DGZ}} or {\cite[p. 137] {G108}}). 


We assume that
the following conditions are satisfied on ${\mathcal U}_{H} \subset M$ 
\begin{eqnarray}
H^{3} &=& \mathrm{tr} (H) \,  H^{2} + \psi \, H + \rho \, g ,
\label{2022.10.01.lla}\\
C \cdot C &=& Q(g, T) ,
\label{2022.10.01.llb}
\end{eqnarray}
where $T$ is a generalized curvature tensor 
and $\psi$ and $\rho$ some functions on ${\mathcal U}_{H}$. 
Now, in view of {\cite[Theorem 4.5] {2020_DGZ}}, we obtain
\begin{eqnarray}
T &=& \left(
\frac{\kappa + 2 \varepsilon \psi }{n-1}
- \frac{\widetilde{\kappa}}{n+1} \right) C
+ \frac{\lambda_{1} }{2}\, g \wedge g \nonumber\\
& &
- \frac{n-3}{ (n-2)^{2} (n-1)}
\left( g \wedge S^{2} + \frac{n-2}{2}\, S \wedge S - \kappa \, g \wedge S
\right) 
\label{2022.10.01.llc}
\end{eqnarray}
on ${\mathcal U}_{H}$, where $\lambda _{1} $ is some function on this set.
Using
(\ref{2022.11.10.aaa}), (\ref{2022.10.01.llb}) and (\ref{2022.10.01.llc})
we get immediately 
\begin{eqnarray*}
T &=& \left(
\frac{\kappa + 2 \varepsilon \psi }{n-1}
- \frac{\widetilde{\kappa}}{n+1} \right) C
- \frac{n-3}{ (n-2)^{2} (n-1)}\, E
+ \frac{\lambda }{2}\, g \wedge g 
\end{eqnarray*}
and
\begin{eqnarray}
C \cdot C &=& 
\left(
\frac{\kappa + 2 \varepsilon \psi }{n-1}
- \frac{\widetilde{\kappa}}{n+1} \right) Q(g,C)
- \frac{n-3}{ (n-2)^{2} (n-1)}\,Q(g,E)
\label{2022.11.01.ll}
\end{eqnarray}
on ${\mathcal U}_{H}$, where $\lambda $ is some function on this set. 
In addition, if we assume that (\ref{4.3.012}) holds on ${\mathcal U}_{H}$
then from  (\ref{2022.10.01.ll}) it follows that
\begin{eqnarray}
\left(
\frac{\kappa + 2 \varepsilon \psi }{n-1} 
- \frac{\widetilde{\kappa}}{n+1} -  L_{C} \right)  C 
&=&  \frac{n-3}{ (n-2)^{2} (n-1)}\,E 
+ \frac{\lambda _{2} }{2}\, g \wedge g 
\label{2023.03.23.aa}
\end{eqnarray}
on ${\mathcal U}_{H}$,
where $\lambda _{2}$ is some function on this set.
We note that (\ref{2023.03.23.aa}), 
by a suitable contraction, yields $\lambda _{2} = 0$,
and in a consequence we obtain
\begin{eqnarray}
\left(
\frac{\kappa + 2 \varepsilon \psi }{n-1} 
- \frac{\widetilde{\kappa}}{n+1} -  L_{C} \right)  C 
&=&  \frac{n-3}{ (n-2)^{2} (n-1)}\, E .
\label{2023.03.23.bb}
\end{eqnarray}
From the above presented considerations it follows immediately 
the following result.
\begin{thm}
Let $M$ be a non-Einstein and non-conformally flat hypersurface 
in $N_{s}^{n+1}(c)$, $n \geq 4$. 
If (\ref{4.3.012}) and (\ref{2022.10.01.lla}) are satisfied on 
${\mathcal U}_{H} \subset M$ then (\ref{2023.03.23.bb}) holds on 
${\mathcal U}_{H}$.
\end{thm}


According to {\cite[Corollary 4.1] {R99}},
if on the subset ${\mathcal U}_{H}$ 
of a hypersurface $M$ in $N_{s}^{n+1}(c)$, $n \geq 4$,
one of the tensors $R \cdot C$, $C \cdot R$ 
or $R \cdot C - R \cdot C$ 
is a linear combination of $R \cdot R$ and of a finite sum 
of tensors of the form $Q(A,T)$, where
$A$ is a symmetric $(0,2)$-tensor 
and $T$ a generalized curvature tensor, then (\ref{2022.10.01.lla}) 
holds on ${\mathcal U}_{H}$.
In particular if one of  the following conditions is satisfied on 
${\mathcal U}_{H} \subset M$: $R \cdot C = Q(g,T_{1})$, 
$C \cdot R = Q(g,T_{2})$ or $R \cdot C - C \cdot R = Q(g,T_{3})$,
where $T_{1}$, $T_{2}$ and $T_{3}$ are generalized curvature tensors, 
then (\ref{2022.10.01.lla}) holds on ${\mathcal U}_{H}$.
Now from Theorems 5.2, 5.3 and 5.4 of \cite{DGPSS},
in view of Proposition 2.2, it follows that
\begin{eqnarray}
\mathrm{Weyl} (T_{1}) &=&
\left(
\frac{\kappa + \varepsilon \psi }{n-1} 
- \frac{(n-1) \widetilde{\kappa}}{n (n+1)} \right)  C 
- \frac{1}{(n-2)(n-1)}\, E ,\label{2023.03.23.cc1}\\
\mathrm{Weyl} (T_{2}) &=&
\left(
\frac{\kappa + 2 \varepsilon \psi }{n-1} 
- \frac{ \widetilde{\kappa}}{n+1} \right)  C 
- \frac{n-3}{(n-2)^{2} (n-1)}\, E ,\label{2023.03.23.cc2}\\
\mathrm{Weyl} (T_{3}) &=&
\left(
 \frac{ \widetilde{\kappa}}{n (n+1)} - \frac{ \varepsilon \psi }{n-1} \right)  C 
- \frac{1}{ (n-2)^{2} (n-1)}\, E .\label{2023.03.23.cc3}
\end{eqnarray}
Thus we have
\begin{thm}
Let $M$ be a non-Einstein and non-conformally flat hypersurface
 in $N_{s}^{n+1}(c)$, $n \geq 4$, satisfying  (\ref{2022.10.01.lla}) 
on ${\mathcal U}_{H} \subset M$,
and let $T_{1}$, $T_{2}$ and $T_{3}$ be generalized curvature tensors 
defined on ${\mathcal U}_{H}$. 
If one of the following conditions: 
$R \cdot C = Q(g,T_{1})$, respectively,
$C \cdot R = Q(g,T_{2})$ and $R \cdot C - C \cdot R = Q(g,T_{3})$, 
is satisfied on ${\mathcal U}_{H}$ 
then (\ref{2023.03.23.cc1}), respectively (\ref{2023.03.23.cc2}) 
and (\ref{2023.03.23.cc3}), holds on ${\mathcal U}_{H}$.
\end{thm}


Finally, we assume that the tensor $H$ satisfies (\ref{2022.10.01.kk}) 
on ${\mathcal U}_{H} \subset M$.
Now (\ref{2022.10.01.ll}), by making use of (\ref{realC6bb}) 
and  (\ref{realRa}), turns into 
\begin{eqnarray}
\frac{n-2}{n-3} \, C \cdot C 
&=& \rho \, Q(g,C)
+  Q \left( S -   \frac{ (n-1) \widetilde{\kappa}}{ n (n+1)}\, g ,  
R - \frac{\widetilde{\kappa}}{2 n (n+1)}\, g \wedge g \right) .
\label{2023.03.24.aa}\\
\rho &=&
\frac{1}{n-3}
\left( \frac{ \kappa }{n-1}  - \frac{  \widetilde{\kappa} }{n+1}  
+ \varepsilon \psi \right) .
\label{2023.03.24.bb}
\end{eqnarray}
From  (\ref{2023.03.24.aa}), by an application of (\ref{abcd07zz}), we obtain
\begin{eqnarray*}
\frac{n-2}{n-3} \, C \cdot C 
&=& \rho \, Q(g,R) + \frac{\rho}{2 (n-2)}\, Q(S, g \wedge g) 
+  Q \left( S -   \frac{ (n-1) \widetilde{\kappa}}{ n (n+1)}\, g ,  
R - \frac{\widetilde{\kappa}}{2 n (n+1)}\, g \wedge g \right) \nonumber\\
&=&
\frac{\rho}{2 (n-2)}\, Q(S,  g \wedge g) 
+ \rho \, 
 Q\left(g,R - \frac{\widetilde{\kappa}}{2 n (n+1)}\, 
 g \wedge g \right) \nonumber\\ 
& & +  Q \left( S -   \frac{ (n-1) \widetilde{ \kappa}}{ n (n+1)}\, g , 
R - \frac{\widetilde{\kappa}}{2 n (n+1)}\, g \wedge g \right)  \nonumber\\
& = &
Q\left( S, \frac{\rho}{2 (n-2)}\, g \wedge g \right) 
 +  Q \left(  S -  \left(  \frac{ (n-1) \widetilde{ \kappa}}{ n (n+1)} 
 - \rho \right) g , R - \frac{\widetilde{\kappa}}{2 n (n+1)}\, 
 g \wedge g \right) 
\nonumber\\
& = &
Q \left( S  -  \left(  \frac{ (n-1) \widetilde{ \kappa}}{ n (n+1)} 
- \rho \right) g ,  \frac{\rho}{2 (n-2)}\, g \wedge g \right) \nonumber\\
& & +  Q \left(  S -  \left(  \frac{ (n-1) \widetilde{ \kappa}}{ n (n+1)} 
- \rho \right) g , R - \frac{\widetilde{\kappa}}{2 n (n+1)}\, g \wedge g \right) 
\nonumber\\
&=&
  Q \left(  S -  \left(  \frac{ (n-1) \widetilde{ \kappa}}{ n (n+1)} 
  - \rho \right) g , 
R - \left( \frac{\widetilde{\kappa}}{ n (n+1)} 
- \frac{\rho}{n-2} \right) \frac{1}{2}\, g \wedge g \right)  .
\end{eqnarray*}
Thus we see that if the tensor $H$ satisfies (\ref{2022.10.01.kk}) 
on ${\mathcal U}_{H} \subset M$ then 
\begin{eqnarray}
\ \ \ \
C \cdot C &=&
 \frac{n-3}{n-2} \,
 Q \left(  S -  \left(  \frac{ (n-1) \widetilde{ \kappa}}{ n (n+1)} 
 - \rho \right) g , 
R - \left( \frac{\widetilde{\kappa}}{ n (n+1)} 
-   \frac{\rho}{n-2} \right) \frac{1}{2}\, g \wedge g \right) 
\label{2023.03.24.cc}
\end{eqnarray}
on  ${\mathcal U}_{H}$, where the function $\rho$ 
is defined by (\ref{2023.03.24.bb}).

In addition, we assume that (\ref{4.3.012}) holds 
on  ${\mathcal U}_{H} \subset M$. Now 
(\ref{2023.03.24.aa}) turns into 
\begin{eqnarray}
& &
\tau \,  Q(g,C) 
- Q\left( S -   \frac{ (n-1) \widetilde{\kappa}}{ n (n+1)}\, g ,  
R - \frac{\widetilde{\kappa}}{2 n (n+1)}\, g \wedge g \right)  \ =\ 0 ,
\label{2023.03.24ee}\\
& &
\tau  \ =\  \rho  - \frac{n-2 }{n-3} \, L_{C} .
\label{2023.03.23.ff}
\end{eqnarray}
From the presented above calculations it follows that 
(\ref{2023.03.24ee}) yields
\begin{eqnarray}
Q\left( S - \left(\frac{ (n-1) \widetilde{\kappa}}{ n (n+1)} 
- \tau  \right) g, R 
- \left( \frac{ \widetilde{\kappa}}{ n (n+1)} 
- \frac{\tau}{n-2} \right) \frac{1}{2}\, g \wedge g \right) &=& 0 .
\label{2023.03.23.gg}
\end{eqnarray}
If 
\begin{eqnarray}
\mathrm{rank} \left(
S - \left(\frac{ (n-1) \widetilde{\kappa}}{ n (n+1)}  
- \tau  \right) g \right) &=& 1
\label{2023.03.23.ww}
\end{eqnarray}
at a point $x \in {\mathcal U}_{H}$ then in view of 
Proposition 2.1
$E = 0$ at $x$, where the tensor $E$ is defined by  (\ref{2022.11.10.aaa}). 
If 
\begin{eqnarray}
\mathrm{rank} \left(
S - \left(\frac{ (n-1) \widetilde{\kappa}}{ n (n+1)} 
- \tau  \right) g \right) & > & 1
\label{2023.03.23.zz}
\end{eqnarray}
at a point $x \in {\mathcal U}_{H}$ then by an application of
{\cite[Proposition 2.4] {DGHHY}}
(or, {\cite[Proposition 2.1] {DGHZ01}})
it follows that the following equation is satisfied at $x$
\begin{eqnarray}
& &
 R - \left(
\frac{ \widetilde{\kappa}}{ n (n+1)} 
- \frac{\tau}{n-2} \right) \frac{1}{2}\, g \wedge g\nonumber\\
&=& \phi 
\left(
S - \left(\frac{ (n-1) \widetilde{\kappa}}{ n (n+1)} 
- \tau  \right) g \right)
\wedge
\left(
S - \left(\frac{ (n-1) \widetilde{\kappa}}{ n (n+1)} 
- \tau  \right) g \right), \ \ \ \phi \in \mathbb{R} .
\label{2023.03.23.yy}
\end{eqnarray}
This by Proposition 4.4 implies 
$(n-2)\, C = \phi \, E$.
Thus we have

\begin{thm}
Let $M$ be a non-Einstein and non-conformally flat hypersurface 
in $N_{s}^{n+1}(c)$, $n \geq 4$. 
\newline
(i) 
If (\ref{2022.10.01.kk})  is satisfied on ${\mathcal U}_{H} \subset M$ then
(\ref{2023.03.24.cc}) holds on ${\mathcal U}_{H}$, where 
the function $\rho$ is defined by (\ref{2023.03.24.bb}) on this set.
\newline
(ii)
If (\ref{4.3.012}) and (\ref{2022.10.01.kk}) are satisfied 
on ${\mathcal U}_{H} \subset M$ then
(\ref{2023.03.23.gg})
holds on ${\mathcal U}_{H}$, where the function $\tau$ 
is defined by (\ref{2023.03.23.ff}) on this set.
Moreover, $\lambda \, C = E$ on  ${\mathcal U}_{H}$, 
where $\lambda$ is some function on this set. 
\end{thm}

\vspace{5mm}

\noindent
\footnotesize{Ryszard Deszcz\\
retired employee of the
Department of Applied Mathematics\\
Wroc\l aw University of Environmental and Life Sciences\\ 
Grunwaldzka 53, 50-357 Wroc\l aw, Poland}\\
{\sf E-mail: Ryszard.Deszcz@upwr.edu.pl} 

\vspace{1mm}

\noindent
\footnotesize{Ma\l gorzata G\l ogowska\\
Department of Applied Mathematics \\
Wroc\l aw University of Environmental and Life Sciences\\ 
Grunwaldzka 53, 50-357 Wroc\l aw, Poland}\\
{\sf E-mail: Malgorzata.Glogowska@upwr.edu.pl}

\vspace{1mm}

\noindent
\footnotesize{Marian Hotlo\'{s}\\
retired employee of the Department of Applied Mathematics\\ 
Wroc{\l}aw University of Science and Technology\\
Wybrze\.{z}e Wyspia\'{n}skiego 27\\ 
50-370 Wroc{\l}aw, Poland}\\ 
{\sf E-mail: Marian.Hotlos@pwr.edu.pl}

\vspace{1mm}

\noindent
\footnotesize{Miroslava Petrovi\'{c}-Torga\v{s}ev\\
Department of Natural and Mathematical Sciences\\
State University of Novi Pazar\\
Vuka Karad\v{z}i\'{c}a bb, Novi Pazar, RS-36300, Serbia}\\
{\sf E-mail: mirapt@kg.ac.rs}

\vspace{1mm}

\noindent
\footnotesize{Georges Zafindratafa\\
Laboratoire de Math\'{e}matiques pour l'Ing\'{e}nieur (LMI)\\
Universit\'{e} Polytechnique Hauts-de-France\\
59313 Va\-len\-cien\-nes cedex 9, France}\\
{\sf E-mail: Georges.Zafindratafa@uphf.fr}

\end{document}